\documentclass[12pt,a4paper,fleqn]{amsart}

\usepackage{geometry}

\usepackage{amsmath}
\usepackage{amsthm, amssymb, mathabx, graphicx}
\usepackage[pdftex,bookmarks=true]{hyperref}
\usepackage{color}
\usepackage{mathtools}
\usepackage{esint}

\newtheorem{theorem}{Theorem}
\newtheorem*{theorem*}{Theorem}

\newtheorem{corollary}{Corollary}
\newtheorem*{corollary*}{Corollary}

\newtheorem*{definition*}{Definition}
\newtheorem{lemma}{Lemma}
\newtheorem*{lemma*}{Lemma}

\newtheorem{proposition}{Proposition}
\newtheorem*{proposition*}{Proposition}

\numberwithin{equation}{section}

\newtheorem{condition}{Condition}
\newtheorem*{condition*}{Condition}
\newtheorem{claim}{Claim}

\def\C{{\mathbb C}}
\def\R{{\mathbb R}}

\def\E{{\mathbb E}}

\def\P{{\mathbb P}}

\def\<{{\langle}}
\def\>{{\rangle}}

\title[Random polynomials with coefficients of polynomial growth]{Real roots of random polynomials with   coefficients of   polynomial growth: a comparison principle and applications}
\dedicatory{In loving memory of Giang T. Ho}
\author{Yen Q. Do}
\address{Department of Mathematics, The University of Virginia, Charlottesville, VA 22904-4137}
\email{yendo@virginia.edu}

\subjclass[2000]{30B20}
\thanks{Y.D. partially supported by NSF grant DMS-1800855.}
\date{}

\begin{document}
\begin{abstract} This paper seeks to further explore the distribution of the real roots of random polynomials with non-centered  coefficients. 
We focus on polynomials where the typical values of the coefficients have power growth and  count the average number of real zeros.  Almost all previous results require coefficients with zero mean, and it is non-trivial to extend these results to  the general case. Our approach is based on a novel comparison principle that reduces the general situation to the mean-zero setting. As applications, we obtain  new results for the Kac polynomials,   hyperbolic random polynomials,  their derivatives, and generalizations of these polynomials. The proof features new logarithmic integrability estimates for random polynomials (both local and global) and fairly sharp estimates for the local number of real zeros. 
\end{abstract}
\maketitle

\tableofcontents

\setlength{\mathindent}{0pt}

\section{Introduction and statement of results}
This paper seeks to  further explore the distribution of the real roots of random algebraic polynomials 
$$p_n(z)= a_0+ a_1z+\dots + a_n z^n, \ \ z\in \C,$$ 
where the coefficients   $a_0,\dots, a_n$ are independent real-valued random variables with finite means  and finite variances.  We are particularly interested in the average number of real roots of such polynomials. This problem has attracted many mathematicians' attention since previous centuries, initially out of theoretical curiosity, but has recently found   applications in  statistical physics and finance  \cite{hkpv2009, sambandham1979, prosen1996, sm2008}.  It was reported in \cite{todhunter2014} that during the $18$th century Waring considered the distribution of the real roots for  random polynomials of low degrees. It however took quite a while until  the first  (but rather crude) estimates for the number of real roots for random polynomials were established, in a  result of   Bloch and Polya at the beginning of the 20th century \cite{bp1932}. Various authors subsequently worked on this problem,  leading to significant developments during 1940s-1970s, with seminal contributions  of Kac~\cite{kac1943}, Littlewood and Offord~\cite{lo1943, lo1945,lo1948}, Ibragimov and Maslova \cite{im1968, im19711, im19712, im19713, maslova1974, maslova1975}, among others.  Recently, there has been a renewed interest in this problem  \cite{ek1995, kz2012, iz2013, dalmao2015, gkz2015,soze20171,soze20172, dv2017, fk2018,nv2019}, in particular Tao and Vu \cite{tv2015} developed a new framework to study the real roots of random polynomials, adapting their methods from random matrix theory.  See also  \cite{nnv2015, dnv2015,dnv2017,nv2017} for  some further development of the methods in \cite{tv2015}.  

Despite the large number of prior studies, only a very  few   are about random polynomials with non-centered coefficients, namely when the coefficients may have \emph{nonzero means}. Furthermore,   these studies often require very restrictive assumptions of algebraic nature on the relationship between the mean, the variance, and the underlying index of   the coefficients.   Ibragimov and Maslova \cite{im19711, im19712} in 1970s considered random  polynomials with iid coefficients of nonzero mean (these are known as Kac polynomials). They showed that the expected number of real roots for the Kac random polynomials  is essentially reduced   to a half if the iid coefficients have a (common) nonzero mean. In \cite{dnv2017},  a joint work with Oanh Nguyen and Van Vu, using different methods we strengthened and generalized  this result  to   random polynomials where the mean and the variance of the coefficient $a_j$ are linearly dependent and furthermore they are  algebraic polynomials of  $j$.  

In this paper, we consider an innovative approach that circumvents the needs for algebraic constraints between the mean and the variance of the coefficients and does not require any algebraic dependence on the underlying index. In particular, this approach  offers some explanation for  the interaction between the mean and the variance of random polynomials. We focus on generalized Kac polynomials, an important class  where  the typical values of the coefficients are  comparable to a fixed power of the underlying index. We will discuss below the technical details of our set up.\footnote{It may be possible that the current approach will be applicable to some other  classes of random functions (such as those studied in \cite{nv2017}), however this will not be explored in this paper and left for further studies.}

For convenience of notation, we  write  $a_j=b_j+c_j\xi_j$
where
$$b_j=\E [a_j] \quad \text{and} \quad |c_j|=\sqrt{Var[a_j]}.$$ 
Note that we do not assume $c_j\ge 0$  and prefer to leave the setup in this generality for the convenience of the proof. Let $\rho\in \R$. For the  typical values  of $|a_j|$ to be comparable to $(1+j)^\rho$, it is natural to assume that $ \E [a_j]= O((1+j)^\rho)$ and  $(Var[a_j])^{1/2}$ is comparable to $(1+j)^\rho$, so that there is a significant range of values for $|a_j|$  about the size of $(1+j)^\rho$. The following condition  essentially describes these assumptions. For technical reasons, below we will need  $\rho>-1/2$.

\begin{condition}\label{cond.polyvar} Assume that for some $\epsilon_0,C_0,N_0>0$ and $\rho>-1/2$ it holds that

(i) $ \E |\xi_j|^{2+\epsilon_0} \le C_0$ for all $0\le j\le n$;  

(ii) $|b_j|, |c_j| \le  C_0(1+j)^{\rho}$ for all $j$;

(iii) $|c_j| \ge \frac 1{C_0} (1+j)^\rho$ for $N_0 \le j\le n-N_0$.

\end{condition}

We note that  $b_j$ and $c_j$ may depend on $n$.  Without loss of generality, we may assume that $0<\epsilon_0\le 1$ throughout the paper. The implicit constants in this paper are allowed to  depend on  the implicit constants in Condition~\ref{cond.polyvar}, which include  $\rho, \epsilon_0,C_0,N_0$.  

We now mention several examples that satisfy Condition~\ref{cond.polyvar}. Via Stirling's formula, it can be seen that the coefficients of hyperbolic random polynomials\footnote{For  discussions about the importance of random hyperbolic polynomials in statistical physics, we refer the reader to the beautiful lecture notes   \cite{hkpv2009}.}
\begin{eqnarray}\label{e.hyperbolic}
p_{\xi,L,n}(z) &=& \xi_0 + \sqrt L \xi_1 z + \dots + \sqrt{\frac{L\dots (L+n-1)}{n!}}\xi_n z^n   
\end{eqnarray}
satisfy the above condition; here  $L>0$ and $\xi_j$'s are independent with unit variance. In particular, if $L=1$ we recover the  Kac random polynomials. 
In fact,  we may generate other examples satisfying Condition~\ref{cond.polyvar} by taking finite linear combinations of   hyperbolic polynomials and their derivatives.  Now, while our approach works with  more general polynomials, even for the polynomials considered in  \cite{dnv2017, im19711, im19712} we are also able to obtain significantly new results.

\subsection{Notational conventions}
Throughout the paper, for any function $q :\R \to \C$ we let $N_{q}$ denote the number of its real roots, and let $N_q(I)$ be the number of  roots inside  $I\subset \R$. Note that these numbers could be $\infty$, but they are never negative. 

By $A\lesssim_{t_1,\dots,} B$ we mean $A=O_{t_1,\dots}(B)$,  in other words there is a finite constant $C$ such that $|A|\le CB$ and the constant $C$ is allowed to depend on the parameters $t_1,\dots$. Sometimes we will simply write $A\lesssim B$ (without mentioning the parameters $t_1,\dots$) when $C$ is an absolute consatnt or if it is clear from the context what $C$ could depend on.    When both $A\lesssim B$ and $B\lesssim A$ hold we will write $A\approx B$, and we use the same convention for $A\approx_{t_1,\dots} B$.

The reciprocal polynomial for a   polynomial $p_n$ of degree $n$ is $p^*_n(z):=z^n p_n(1/z)$.

\subsection{Statement of results}

To study $N_{p_n}$, we write
$$p_n(z) = m_n(z)+r_n(z)$$
where $m_n(z)=\E p_n(z)$ is a deterministic polynomial and $r_n = p_n-m_n$ is a random polynomial with zero mean.  Our heuristics is the following idea: locally, between $m_n$ and $r_n$, the dominant component  will dictate the behavior of $p_n$ and hence will have a stronger influence on  the number of real zeros of $p_n$. 

Our main result, Theorem~\ref{t.general} is an estimate for the number of real roots of $p_n$ inside an arbitrary interval, demonstrating  the following comparison principle:  

(i) if $m_n$ dominates $r_n$    then on average there are very few real roots for $p_n$, as $|m_n|$ is typically bigger than $|r_n|$. 

(ii) if $m_n$ is dominated by $r_n$  then on average the  number of real roots of $p_n$  is the same as the  number of real roots of $r_n$ plus a bounded term.

In the statement of Theorem~\ref{t.general}, we will be more precise about the meaning of ``dominated'' and ``dominates''. Here we make some preliminary remarks. First, since $r_n$ is random with zero mean, it makes sense to use the standard deviation $(Var[r_n])^{1/2}$  as an indicator for the size of $r_n$, and this heuristics is also used for derivatives of $r_n$.  For  $t\ge 1$, to compare $m_n$ and $r_n$  it turns out to be more convenient to work with the reciprocal polynomials $m^*_n$ and $r^*_n$.

In the following,  we say that $J$ is an enlargement  for  $I=(a,b)$ if it is obtained by extending $I$ to the left and to the right a little bit: generally speaking this means there is an absolute constant $c>0$ such that   the added length to the right is bounded below by $c(\Big|1-|b|\Big|+\frac 1 n)$ and the added length to the left is bounded below by $c(\Big|1-|a|\Big|+\frac 1 n)$.   

There are special cases when the enlargement requirement could be made less stringent (without affecting our main results below): if $|1-|b||$ is bounded below by any positive absolute constant then there is no need to extend $I$ to the right and we may use $b$ as the right endpoint for $J$,  and similarly if $|1-|a||$ is bounded below by any positive absolute constant then we may take $a$ as the left endpoint for $J$.   These improvements are made possible with the aid of Lemma~\ref{l.im}.

We note that the above notion of enlargement can also be similarly defined for  half open/half closed/closed/infinite intervals.  In all cases, the following will be true: if $J$ is an enlargement of $I$ then it also qualifies as an enlargement of any subintervals of $I$.


  
 \begin{theorem}[Comparison principle]\label{t.general} There is a  constant $0<C  <\infty$ such that the following holds. Assume that the coefficients of $p_n$ satisfy Condition~\ref{cond.polyvar} and are real valued. Let $I\subset \R$ be an interval whose endpoints may depend on $n$ and assume that  $J$ is an enlargement of $I$.

Let $m_n^*(t)=t^n m_n(\frac 1 t)$ and $r^*_n(t)=t^n r_n(\frac 1 t)$ for $t\ne 0$.

(1) Assume that
\begin{itemize}
\item  if $t\in   J\cap [-1,1]$ then $\quad |m_n(t)|   >   C |\log(1-|t|+\frac 1 n)|^{1/2} \sqrt{Var[r_n(t)]}$,
\item  if $t\in   J\setminus [-1,1]$ then $|m^*_n(\frac 1 t)|   >   C |\log(1-\frac 1{|t|}+\frac 1 n)|^{1/2} \sqrt{Var[r^*_n(\frac 1 t)]}$.
\end{itemize}

Then $\E N_{p_n}(I) = O(1).$

 (2) Let $\phi: [0,1]\to [0,1]$ such that    $\displaystyle \int_{1/n}^c \frac{\phi(t)}tdt = O(1)$  for some $c>0$.
 
Assume that for each $k=0,1$ we have the uniform estimates:
\begin{itemize}
\item if  $t\in J \cap[-1,1]$ then $\quad |m^{(k)}_n(t)|  \lesssim  \phi(1-|t|+\frac 1 n)\sqrt{Var[r^{(k)}_n(t)]}$, 
\item if  $t\in J\setminus [-1,1]$ then $\quad |{(m^*_n)}^{(k)}(\frac 1 t)|   \lesssim  \phi(1-\frac 1{|t|}+\frac 1n) \sqrt{Var[{(r_n^*)}^{(k)}(\frac 1 t)]}$,
\end{itemize}
and for $k=2$ the weaker  estimates without  $\phi$ also hold on  $J \cap[-1,1]$ and $J\setminus [-1,1]$.

Then $\E N_{p_n}(I) = \E N_{r_n}(I) + O(1).$

\end{theorem}
  
We note that Theorem~\ref{t.general} is more useful for intervals near $\pm 1$, since under Condition~\ref{cond.polyvar} it can be shown (using a standard argument of Ibragimov and Maslova) that $\E N_{p_n}(I)=O(1)$ if $I$ is bounded away from $\pm 1$ 
(see Lemma~\ref{l.im}).  

In Theorem~\ref{t.general}, for technical reasons we need to assume that the domination  relationship (between $m_n$ and $r_n$) is effective on an enlargement $J$ of $I$, however  if $p_n$ is a Gaussian random polynomial then the conclusions hold  with $J=I$ and some of the conditions could be weakened, see Section~\ref{s.gaussian}. The proof of the Gaussian case in Section~\ref{s.gaussian} will also shed more light on the motivation for the  assumptions on $m_n$ and $r_n$ in the statement of Theorem~\ref{t.general}.  One of the main technical ingredients in our proof is a new result about universality for the correlation of the roots of $p_n$,  see Section~\ref{s.corr}.

Using Theorem~\ref{t.general}, we could derive new results about   the real roots of  non-centered random polynomials (with coefficients of power growth) from analogous results for centered random polynomials, which in turn were studied extensively in \cite{dnv2017}. Below, we summarize several sample results that can be obtained in this direction  (although this list is by no means  comprehensive).\footnote{A more thorough discussion  about possible applications is included in Section~\ref{s.applications}, where these sample results will be derived from Theorem~\ref{t.general}.}  The sample results  will further demonstrate  the following observation from \cite{dnv2017}: we may extract  asymptotic estimates for the number of real roots of a random polynomial from asymptotic information about its coefficients. This phenomenon was first observed in \cite{dnv2017}  for   random polynomials with centered coefficients of polynomial growth.

Below, following \cite{dnv2017}, we define a generalized polynomial of $j\in \mathbb Z_+$ to be a finite linear combination of hyperbolic coefficients $h_L(j):=\frac{L(L+1)\dots (L+j-1)}{j!}$, $L>0$. Its degree is defined to be  $L_{max}-1$, where $L_{max}$ is the biggest $L$ in the combination. If we requires $L$ to be integer then this notion is the same as the classical notion of polynomials. Note that (via Stirling's formula)  a generalized polynomial of degree $\delta$ is asymptotically comparable to $j^{\delta}$.

Our first sample result is about random hyperbolic polynomials \eqref{e.hyperbolic}.

\begin{theorem}\label{t.hyperbolic}  Let $p_n$ be the hyperbolic random polynomial $p_{\xi, L,n}$ given by \eqref{e.hyperbolic} where $\xi_j$ are independent with a common nonzero mean and variance $1$ and uniformly bounded $(2+\epsilon)$ moments  for some $\epsilon>0$. 
\begin{align*}
\text{Then}&& \E N_{p_n} = \frac{(1+\sqrt {L})\log n}{2\pi}+O(1),
\end{align*}
\begin{align*}
\text{and for any $k\ge 1$  we have}&& \E N_{p^{(k)}_n} = \frac{(1+\sqrt {L+2k})\log n}{2\pi}+o(\log n).
\end{align*}
\end{theorem}

Theorem~\ref{t.hyperbolic} is a special case of the following  more general result.

\begin{theorem}\label{t.hyperbolic-gen}  Assume that the coefficients of $p_n$ satisfy Condition~\ref{cond.polyvar}.  Assume furthermore that there are   $\rho_1<\rho-1/2 < \rho_2$ such that   $|b_j| \gtrsim j^{\rho_2} +O(1)$ and
$$|b_{j+1} - b_j| = O((j+1)^{\rho_1}).$$ 
\begin{align*}
\text{Then for any $C>0$ we have}&&
\E N_{p_n} =  \E N_{r_n}(1-\frac 1 C, 1+\frac 1 C)+ O(1),
\end{align*}
in particular $\E N_{p_n}$ grows like $\log n$ as $n\to \infty$.
Furthermore, if  for some $C$ we have $c_j = (C+o(1))j^\rho$ as $j\to\infty$ then
$$\E N_{p_n} =  \frac{1+\sqrt{2\rho+1}}{2\pi} \log n  + o(\log n).$$
In particular, if $c_j^2$ is a generalized polynomial  of $j$ then
$$\E N_{p_n} =  \frac{1+\sqrt{2\rho+1}}{2\pi} \log n  + O(1).$$
\end{theorem}
Theorem~\ref{t.hyperbolic} may be derived from Theorem~\ref{t.hyperbolic-gen} as follows. Letting $\rho=(L-1)/2$,  we note that for the set up of Theorem~\ref{t.hyperbolic} we will have $b_j=c_j \mu$ for some $\mu\ne 0$, and by Stirling's formula $c_j=\sqrt{\frac{L(L+1)\dots (L+j-1)}{j!}} = (C_L+o(1))(1+j)^{\rho}$. On the other hand,
$$|b_{j+1}-b_j| = |\mu c_j| |\sqrt{(L+j)/(j+1)} - 1|  = O((j+1)^{-1}c_j)=O((j+1)^{\rho-1}).$$
Using Theorem~\ref{t.hyperbolic-gen}, it follows that $\E N_{p_n} = \E N_{r_n}(1-\frac 1 C, 1+\frac 1C)$, and thus using \cite{dnv2017} we obtain the desired conclusions. We may argue similarly to get the desired  asymptotics for $\E N_{p^{(k)}_n}$.

Below is a class  of random polynomials where the deterministic component $m_n$ is  dominated by the random component $r_n$.

\begin{theorem}\label{t.asymplogs} Assume  Condition~\ref{cond.polyvar} and  assume  that for some $\rho'<\rho-1/2$ we have $|b_j|=O((1+j)^{\rho'})$. Then there are finite positive constants $C_1$ and $C_2$ such that
$$C_1\log n + O(1) \le \E N_{p_n} \le C_2 \log n + O(1).$$
Furthermore if  for some $C$ we have $c_j = (C+o(1))j^\rho$ as $j\to\infty$ then we could take $C_1,C_2$ to be $\frac{1+\sqrt{2\rho+1}}{\pi} + o(1)$.  In particular, if $c_j^2$ is a generalized polynomial  of $j$ then we could let $C_1,C_2=\frac{1+\sqrt{2\rho+1}}{\pi}$.
\end{theorem}

Finally, we mention a simple class of  random polynomials where $m_n$ dominates $r_n$, leading to very few real zeros for the random polynomial.

\begin{theorem}\label{t.asymplogl}Assume Condition~\ref{cond.polyvar}. Suppose furthermore that for some $\rho'\in (\rho-\frac 1 2, \rho]$ and some $\rho'' <\rho'$ the following holds: for odd $j$ we have $b_j = O((1+j)^{\rho''})$ and for even $j$  we have  $b_{j}\gtrsim (1+j)^{\rho'} - O(1)$. Then 
$$\E N_{p_n} = O(1).$$
Furthermore, the above estimate holds true if we interchange the role of odd and even $j$'s in the above assumptions.
\end{theorem}

\subsection{Outline of the paper}
In the next section, we discuss the applications of Theorem~\ref{t.general} and the proof for the sample results mentioned above. In the rest of the paper, we prove Theorem~\ref{t.general}. Our proof of Theorem~\ref{t.general} uses universality estimates for the  correlation functions of the real roots of $p_n$, see Section~\ref{s.corr}. Using these estimates, we could  reduce the proof of Theorem~\ref{t.general} to the Gaussian setting. The Gaussian case of Theorem~\ref{t.general} will be examined using the Kac-Rice formula, see Section~\ref{s.gaussian}.

\section{Sample applications of the comparison principle}\label{s.applications}
In this section, we discuss several applications of Theorem~\ref{t.general} and present the proofs for Theorem~\ref{t.hyperbolic-gen}, Theorem~\ref{t.asymplogs}, and Theorem~\ref{t.asymplogl}. We will use the following basic computation about power series.
\begin{lemma}\label{l.elementary} For any $\alpha>-1$ and $\beta>-1$ and any  $c>0$ and $C>1$ the following holds:

(i) If $\frac 1 C\le t \le 1-\frac c n$ then $\sum_{j=1}^n (n+1-j)^\beta j^\alpha t^j  \ \ \approx_{\alpha, \beta, c,C} \ \  n^{\beta} (1-t)^{-\alpha-1}$.

(ii) If $|1-t|\le c/n$ then  $\sum_{j=1}^n (n+1-j)^\beta j^\alpha t^j  \ \ \approx_{\alpha,\beta,c,C} \ \ n^{\alpha+\beta+1}$.
\end{lemma}

\proof[Proof of Lemma~\ref{l.elementary}] Note that if $1-c/n \le t \le 1+c/n$ then $1, t,\dots, t^n$ are all comparable to $1$, therefore $\sum_{j=1}^n (n+ 1- j)^\beta j^\alpha t^j \approx \sum_{j=1}^n (n+1-j)^\beta j^\alpha \approx n^{\alpha+\beta+1}$.  Here, to see the last estimate we may split the sum into $1\le j\le n/2$ and $n/2<j\le n$, and use the fact that for the first range $n+1-j\approx n$ and for the second range $j\approx n$. This proves part (ii), and furthermore in part (i) we may assume that $1/C\le t\le 1-c/n$ where $c$ is sufficiently large. We now discuss the proof of part (i) under this assumption.

We consider first the case $\beta=0$. By Taylor's theorem, we have $(1-t)^{-\alpha-1} =1+(\alpha+1)t+\dots +  \frac{(\alpha+1)\dots (\alpha+n)}{n!} t^n + E_n(t)$, where the error term $E_n(t)$ is nonnegative. Now, note that  $(\alpha+1)\dots (\alpha+j)/j! \approx j^{\alpha}$, therefore
$$\sum_{j=1}^{n}  j^\alpha t^j \lesssim  (1-t)^{-\alpha-1}.$$
For the other direction of the estimate, it suffices to establish that the error term  $E_n(t)$ is smaller than fraction of $(1-t)^{-\alpha-1}$ when $c$ is sufficiently large.  Here we use the  Lagrange form of the error term, which says that  for some  $s \in (0,t)$ we have
\begin{eqnarray*}
E_n(t) &=&(1-s)^{-\alpha-n-2} \frac{(\alpha+1)\dots (\alpha+n+1)}{(n+1)!} (t-s)^{n+1} \\
&\lesssim& (1-s)^{-\alpha-n-2}   (n+1)^{\alpha} (t-s)^{n+1} \\
&=& (1-t)^{-\alpha-1}   (1-\frac{1-t}{1-s})^{n+1} (\frac{1-t}{1-s})^{\alpha+1} (n+1)^\alpha
\end{eqnarray*}
The desired estimate then follows from the fact that $(1-v)^{n} v^\alpha n^\alpha$ is a decreasing function for $v\in [\alpha/n,1]$, and 
$$(1-c/n)^n (c/n)^\alpha n^\alpha \le e^{-c} c^\alpha $$
and  $e^{-c}c^\alpha$ could be made arbitrarily small by choosing $c$ sufficiently large. 

We now consider the general situation. We have
\begin{eqnarray*}
\sum_{j=1}^{n/2} (n+1-j)^\beta j^\alpha t^j  
&\approx& n^{\beta} \sum_{j=1}^{n/2}  j^\alpha t^j  \quad \approx \quad n^{\beta} (1-t)^{-(\alpha+1)}.
\end{eqnarray*}
Thus it remains to show that the remaining  summation over $n/2<j\le n$ is $O(n^{\beta} (1-t)^{-(\alpha+1)})$ (note that this summation is nonnegative). For these $j$'s we note that $j$ is comparable to $n$. Since $\beta>-1$ we may choose  $1<p<\infty$  depending on $\beta$ such that $\beta p>-1$. Let $q=p/(p-1)$ be its conjugate exponent. Then using H\"older's inequality we have
\begin{eqnarray*}
\sum_{j=n/2}^n (n+1-j)^\beta j^\alpha t^j   
&\lesssim&   (\sum_{j=n/2}^n (n+1-j)^{p\beta})^{1/p} (\sum_{j=n/2}^n j^{q\alpha} t^{qj})^{1/q}\\
&\lesssim&  n^{\beta+1/p} (\sum_{j=n/2}^n j^{q\alpha}t^{qj})^{1/q}\\
&\approx&  n^{\beta} (\sum_{j=n/2}^n j^{q(\alpha+ 1)-1}t^{qj})^{1/q}\\
&\lesssim& n^{\beta} ((1-t)^{-q(\alpha+1)})^{1/q} \quad = \quad n^{\beta} (1-t)^{-(\alpha+1)}.
\end{eqnarray*}
This completes the proof of Lemma~\ref{l.elementary}.
\endproof

Let $C>0$ be a sufficiently large constant  and let $A_C=\{z\in \R:  ||z|-1| > 1/C\}$. In the applications of Theorem~\ref{t.general}, we will need the following estimate.
\begin{lemma}\label{l.im} For any $C>0$ we have
$\E N_{p_n}(A_C) =O_C(1)$.
\end{lemma}
We include a proof  of Lemma~\ref{l.im} using an argument of Ibragimov--Maslova \cite{im19713} (see also \cite{dnv2017} where a simpler version of Lemma~\ref{l.im} was proved). We'll need the following estimate, which will also be used later in the proof of Theorem~\ref{t.general}.
\begin{lemma}\label{l.unif-smallball} For any $\delta_0<1$ there is  $p_0\in (0,1)$ such that for  any $\alpha$ we have $\max_j \P(|\xi_j-\alpha| \le \delta_0) \le 1-p_0$.
\end{lemma}

\proof[Proof of Lemma~\ref{l.unif-smallball}] 
Let $\delta_0<1$ and let $0 \le j \le n$. 

We first consider $|\alpha| > 3$. Without loss of generality assume $\alpha>3$, the case $\alpha<-3$ is can be treated similarly. Then
\begin{eqnarray*}
\P(|\xi_j- \alpha| \le \delta_0) 
&\le& \P(\xi_j \ge \alpha - \delta_0) \\
&\le& (\alpha-\delta_0)^{-2} \E |\xi_j|^2 \le 1/4.
\end{eqnarray*}
Thus we may take any $p_0 \le 3/4$ for $|\alpha|>3$.

We now consider $|\alpha| \le 3$. Then $\E |\xi_j-\alpha|^{2+\epsilon_0} = O_{C_0,\epsilon_0}(1)$. Therefore,
\begin{eqnarray*}
\E |\xi_j-\alpha|^2  &\le&  \delta_0^2 \P(|\xi_j-\alpha| \le \delta_0) + \E [|\xi_j|^2 1_{|\xi_j -\alpha| > \delta_0}]  \\
&\le&   \delta_0^2 \P(|\xi_j-\alpha| \le \delta_0) + (\E |\xi_j-\alpha|^{2+\epsilon_0})^{\frac 2{2+\epsilon_0}} \Big(\P(|\xi_j-\alpha| >\delta_0)\Big)^{\frac{\epsilon_0}{2+\epsilon_0}}.
\end{eqnarray*}
Let $x=\P(|\xi_j-\alpha| > \delta_0) \ge 0$. Since $\E |\xi_j-\alpha|^2 = 1+|\alpha|^2 \ge 1$, we obtain
$$0<1-\delta_0^2 \le C_1 x^{\frac{\epsilon_0}{2+\epsilon_0}} - \delta_0^2 x$$
for some $C_1=C_1(C_0,\epsilon_0)$ where $C_0$ and $\epsilon_0$ are as in Condition~\ref{cond.polyvar}.
Thus by examining the function $C_1x^{\epsilon_0/(2+\epsilon_0)} - \delta_0^2 x$ of $x$,  it is follows that  there is some $p_0=p_0(\delta_0,C_1,\epsilon_0) \in (0,1)$ such that any $x\in [0,1]$ that satisfies the above inequality must be inside $[p_0,\infty)$.  Consequently $\P(|\xi_j-\alpha|\le \delta_0)\ge p_0$, as desired.
\endproof

\proof[Proof of Lemma~\ref{l.im}]  It suffices to show that for $r_1<1$ we have $N_{p_n}(-r_1,r_1) = O_{r_1}(1)$ and $N_{p^*_n}(-r_1, r_1)= O_{r_1}(1)$. We will show in detail the first estimate, and comment on the needed changes for the second estimate.

Take any $r_2\in (r_1,1)$. 
 Let $\delta_0, p_0$ be as in Lemma~\ref{l.unif-smallball}.  
From Condition~\ref{cond.polyvar}, let $j_0$ be such that $c_j \approx (1+j)^\rho$ for $j_0\le j\le n-j_0$. Define  
$$A_k:=\{|\xi_j+\frac{b_j}{c_j}| \le \delta_0, \ \  \forall j_0\le j \le k-1\} \cap \{|\xi_k + \frac{b_k}{c_k}| > \delta_0\}$$ 
for each $j_0\le k\le n-j_0$, and define $A_{n-j_0+1} =\{|\xi_j + \frac{b_j}{c_j}| \le \delta_0, \forall j_0\le j \le n-j_0\}$.

For $k=n - j_0+1$ it is clear that we have $\E [1_{A_k} N_{p_n}(-r_1,r_1)] \le n p_0^{n-2j_0} = O(1)$. 

 For $j_0\le k\le n-j_0$, we have $\P(A_k)\le p_0^{k-j_0}$, thus it suffices to show that
\begin{eqnarray*}
   \E [1_{A_k} N_{p_n}(-r_1,r_1)] \lesssim  k(\log k)  \P(A_k),
\end{eqnarray*}
On the event $A_k$,    we have  $|p_n^{(k)}(0)| =  k! |b_k+ c_k \xi_k|\gtrsim  k!  |c_k| \ \gtrsim\  (k+1)^{\rho}$, thus using Jensen's formula we have
\begin{eqnarray*}
N_{p_n}(-r_1,r_1) 
&\le&   1+k+N_{p^{(k)}_n}(-r_1, r_1)  \ \le \  1+ k + O\Big(\sup_{|z|=r_2} \log |p^{(k)}_n(z)|\Big).
\end{eqnarray*}
Let $n_0$ be an integer larger than $\max(0,\rho)$. Using convexity and Jensen's inequality, we have
\begin{eqnarray*}
\frac 1 {\P(A_k)} \E [1_{A_k} N_{p_n}(-r_1,r_1)]
&\lesssim&  1+k + \log\Big(\P(A_k)^{-1} \E [\sup_{|z|=r_2}  |p_n^{(k)}(z)|] \Big)\\
&\lesssim& 1+k + \log\Big(\sum_{i = 0}^{n-k}  (i+1)\dots (i+k+n_0)  r_2^{i} \Big)   \\
&\le& 1+k +\log(\frac{(k+n_0)!}{(1-r_2)^{k+1+n_0}})  \\
&\lesssim& 1+k \log k.
\end{eqnarray*}

To estimate  $\E N_{p^*_n}(-r_1,r_1)$, we proceed similarly, and  the following estimate will be needed:
$$\E  \sup_{|z|=r_2} |{p^*}^{(k)}_n(z)|   \ \lesssim_{r_2,\rho} \ (n+1-k)^{\rho}  ((2k+1)!)^{1/2} (1-r_2^2)^{-(k+1)},$$
where  $r_2\in (r_1,1)$. To see this estimate, we note that
$$\E  \sup_{|z|=r_2}  |{p^*}^{(k)}_n(z)|   \le \sum_{i>(n-k)/2} (n+1-k-i)^\rho (i+1)\dots (i+k) r_2^{i},$$
then we split the sum into $i\le (n-k)/2$ and $i>(n-k)/2$ and argue as in the proof of Lemma~\ref{l.elementary}. The treatment of $i\le (n-k)/2$ is entirely similar as before, but for $i>(n-k)/2$  we actually need to be more careful (than the proof of Lemma~\ref{l.elementary}) about the dependence on $k$ of the implicit constant. We include the details below. By Cauchy--Schwartz we have
\begin{eqnarray*} 
\sum_{i>(n-k)/2} &\le& (\sum_{i>(n-k)/2} (n+1-k-i)^{2\rho})^{1/2} (\sum_{i>(n-k)/2} (i+1)^2\dots (i+k)^2 r_2^{2i})^{1/2} \\
&\lesssim& (n+1-k)^{\rho+1/2} (\sum_{i>(n-k)/2} (i+1)^2\dots (i+k)^2 r_2^{2i})^{1/2} \\
&\lesssim& (n+1-k)^{\rho} (\sum_{i>(n-k)/2} (i+1)\dots (i+2k+1)  r_2^{2i})^{1/2} \\
&\lesssim& (n+1-k)^\rho ((2k+1)!)^{1/2} (1-r_2^2)^{-(k+1)}.
\end{eqnarray*}
\endproof

We now divide the discussion of the applications of Theorem~\ref{t.general} into three sections, corresponding to whether $m_n$ is always small, or always large, or mixed large/small, in comparison  to $r_n$.

\subsubsection{Small mean}
Here the mean $m_n$ will be  completely dominated by $r_n$. We first state a   corollary of Theorem~\ref{t.general} in this direction, before proving Theorem~\ref{t.asymplogs}.

\begin{corollary}\label{c.smallmean} Let $\phi:[0,1] \to [0,1]$ such that $\int_{1/n}^c \frac{\phi(t)}t dt =O(1)$ for some $c>0$. Assume  Condition~\ref{cond.polyvar} and assume that there is a constant $C>1$  such that for $1/C\le |t|\le 1$ and $0\le k\le 1$ we have 
\begin{eqnarray}\label{e.smallmean}
|m_n^{(k)}(t)| &\lesssim& \phi(1 + \frac 1 n - |t|) (1  +\frac 1 n - |t|)^{-(\rho+k+\frac 1 2)},\\
\nonumber |{m^*}^{(k)}_n(t)| &\lesssim&  n^{\rho} \phi(1  +\frac 1 n - |t|)(1  +\frac 1 n - |t|^{-(k+\frac 1 2)},
\end{eqnarray}
and  assume that the weaker estimates without $\phi$ also hold true for $k=2$.
Then there are finite positive constants $C_1$ and $C_2$ such that
$$C_1\log n + O(1) \le \E N_{p_n} \le C_2 \log n + O(1).$$
Furthermore if  for some $C$ we have $c_j = (C+o(1))j^\rho$ as $j\to\infty$ then we could take $C_1,C_2$ to be $\frac{1+\sqrt{2\rho+1}}{\pi} + o(1)$.  In particular, if $c_j^2$ is a generalized polynomial of $j$ then we could let $C_1,C_2=\frac{1+\sqrt{2\rho+1}}{\pi}$.
\end{corollary}
Thanks to \cite{dnv2017}, the zero-mean  case (i.e. $b_j = 0$ for all $j$) of the above corollary already holds true. Thus, using Lemma~\ref{l.im}  and Theorem~\ref{t.general}, Corollary~\ref{c.smallmean} is a simple consequence of the following estimates
\begin{eqnarray}\label{e.var_rn}
\sqrt{Var[r_n^{(k)}(t)]}  &\approx&  (1+\frac 1 n -|t|)^{-(\rho+k+\frac 1 2)}, \\
\nonumber \sqrt{Var[{r^*}^{(k)}_n(t)]}&\approx&  n^{\rho}(1+\frac 1 n -|t|)^{-(k+\frac 12)},
\end{eqnarray}
which  follows from elementary computations (see  Lemma~\ref{l.elementary} for details).

We now prove Theorem~\ref{t.asymplogs}. Since $\rho>-1/2$, we may assume without loss of generality that $\rho'>-1$. Using Lemma~\ref{l.elementary}, for $|t|\le 1$ we then have  
$$|m^{(k)}_n(t)|\lesssim (1+\frac 1 n - |t|)^{-(\rho'+k+1)}, \quad  |{m^*}^{(k)}_n(t)| \lesssim \frac{n^{\rho'}}{(1-|t|+\frac 1n)^{k+1}},$$
which clearly implies \eqref{e.smallmean}. Thus Theorem~\ref{t.asymplogs} follows from Corollary~\ref{c.smallmean}.

\subsubsection{Large mean}
Here near $\pm 1$ the mean $m_n$ will always dominate $r_n$. As before, we state a corollary of Theorem~\ref{t.general} before proving Theorem~\ref{t.asymplogl}.

\begin{corollary}\label{c.largemean} Let $\varphi: (0,\infty)\to [0,\infty)$ be such that $\varphi(t)\to \infty$ as $t\to 1/n$. Assume  Condition~\ref{cond.polyvar} and assume that there is a constant $C>1$  with the following properties: for $1- \frac 1 C \le |t|\le 1$  we have 
\begin{eqnarray}\label{e.smallmean}
|m_n(t)| &\gtrsim& \varphi(1+\frac 1n-|t|) (1+\frac 1 n - |t|)^{-(\rho+\frac 12)},\\
|m^*_n(t)| &\gtrsim&  n^\rho \varphi(1+\frac 1n-|t|)   (1+\frac 1 n- |t|)^{-\frac 1 2}.
\end{eqnarray}
\begin{align*}
\text{Then} && \E N_{p_n} = O(1).
\end{align*}
\end{corollary}
This corollary follows immediately from \eqref{e.var_rn} and Theorem~\ref{t.general} and Lemma~\ref{l.im}. We now apply  this corollary with $\varphi(t)=t^{-\epsilon}$ to prove Theorem~\ref{t.asymplogl}. By splitting $m_n=m_{n,odd} + m_{n,even}$ and using Lemma~\ref{l.elementary} to treat each of them individually, we obtain (for $1-1/C\le |t|\le 1$)
$$m_n(t) \approx    (1+1/n-|t|)^{-(\rho+1+\epsilon)}, \quad m^*_n(t) \approx    \frac{(n+1)^{\rho'}}{1+1/n-|t|} \gtrsim \frac{n^\rho} {(1+\frac 1 n- |t|)^{\epsilon+\frac 1 2}}$$
where $\epsilon=\rho'+1/2-\rho>0$. Thus Theorem~\ref{t.asymplogl} follows from Corollary~\ref{c.largemean}.

\subsubsection{Mixed case}
Here we consider the mixed situation, where $m_n$ is dominated by $r_n$ on a part of the real line and   dominates $r_n$ elsewhere. In our opinion this is the most interesting case. Here we describe  a simple scenario, which applies to     random Kac polynomials with non-centered coefficients  (considered in \cite{im19712}) as well as linear combination of derivatives of a random Kac polynomial (considered in \cite{dnv2017}), and also hyperbolic random polynomials  with non-centered coefficients (Theorem~\ref{t.hyperbolic} of the current paper). In this scenario, $m_n$ is dominated by $r_n$ near $-1$   while being the dominant component near $1$. (Note that due to symmetry we could also state a symmetric version where the roles of $1$ and $-1$ are interchanged.) 

\begin{corollary}\label{c.mixedmean} Let $\varphi: (0,\infty)\to [0,\infty)$ be such that $\varphi(t)\to \infty$ as $t\to 1/n$. Let $\phi:[0,1]\to [0,1]$ such that $\int_{1/n}^c \frac{\phi(t)}t dt = O(1)$ for some $c>0$.  Assume  Condition~\ref{cond.polyvar} and assume that there is a constant $C>1$  with the following properties: 

(i) for $1- \frac 1 C \le t\le 1$  we have 
\begin{eqnarray}\label{e.mixedmean1}
|m_n(t)| &\gtrsim& \varphi(1+\frac 1n-t) (1+\frac 1 n - t)^{-(\rho+\frac 12)},\\
\nonumber |m^*_n(t)| &\gtrsim&  n^\rho \varphi(1+\frac 1n-t)   (1+\frac 1 n- t)^{-\frac 1 2},
\end{eqnarray}

(ii) for $-1\le t\le -1+\frac 1 C$ and for each $k=0,1$ we have
\begin{eqnarray}\label{e.mixedmean-1}
|m_n^{(k)}(t)| &\lesssim& \phi(1 + \frac 1 n + t) (1  +\frac 1 n +t)^{-(\rho+k+\frac 1 2)},\\
\nonumber |{m^*}^{(k)}_n(t)| &\lesssim&  n^{\rho} \phi(1  +\frac 1 n  + t)(1  +\frac 1 n +t)^{-(k+\frac 1 2)}.
\end{eqnarray}
and the  weaker estimates without $\phi$ also hold true for $k=2$. Then
$$\E N_{p_n}  = \E N_{r_n}(1-1/C,1+1/C) + O(1)$$
and in particular there are constants $C_1,C_2>0$ such that 
$$C_1\log n + O(1) \le \E N_{p_n} \le C_2 \log n + O(1).$$
Furthermore if  for some $C$ we have $c_j = (C+o(1))j^\rho$ as $j\to\infty$ then we could take $C_1,C_2$ to be $\frac{1+\sqrt{2\rho+1}}{2\pi} + o(1)$.  In particular, if $c_j^2$ is a generalized polynomial of $j$ then we could take $C_1 = C_2=\frac{1+\sqrt{2\rho+1}}{2\pi}$.
\end{corollary}
Now, it was shown in \cite{dnv2017} that  $\E N_{r_n}(1-1/C,1+1/C)$ grows like $\log n$, and furthermore if $c_j=(C+o(1))j^\rho$ then $\E N_{r_n}(1-1/C,1+1/C)=\frac{1+\sqrt{2\rho+1}}{2\pi}\log n + o(\log n)$, and the error term could also be improved to $O(1)$ if $c_j^2$ is a generalized polynomial of $j$.  Thus, Corollary~\ref{c.mixedmean} is an immediate consequence of Theorem~\ref{t.general} and \eqref{e.var_rn}.

We now discuss the proof of  Theorem~\ref{t.hyperbolic-gen}. From the given assumption it follows that $b_j$ are of the same sign for $j\gtrsim 1$, so without loss of generality we may assume that $b_j>0$ for $j\gtrsim 1$. Now, using $b_j\gtrsim j^{\rho_2}$ and $\rho_2>\rho-1/2$ one may show that   $m_n(t)$ dominates $r_n(t)$ near $1$. Indeed, by elementary computations (see Lemma~\ref{l.elementary}), for $t\in [1-1/C,1]$ we have
\begin{eqnarray*}
m_n(t) &\gtrsim& (1+\frac 1 n-t)^{-(\rho_2+1)} \gtrsim  (1+\frac 1 n - t)^{\epsilon} \sqrt{Var[r_n(t)]},\\
m^*_n(t) &\gtrsim& n^{\rho_2} (1+\frac 1 n-t)^{-1} \gtrsim (1+\frac 1 n-t)^{\epsilon} \sqrt{Var[r^*_n(t)]}.\end{eqnarray*}
We now show that $m_n$ is dominated by $r_n$ near $-1$. To see this, let $k\ge 0$ and we use discrete integration by parts to write 
$$(k!)^{-1}m^{(k)}_n(t) =  b_k(1+t+\dots+t^{n-k}) +  \sum_{j=k+1}^{n} \Big({j\choose k}b_{j}- {j-1 \choose k}b_{j-1}\Big)(t^{j-k}+\dots+t^{n-k})$$
and uniformly over $j_1\le j_2$ we have $t^{j_1}+\dots + t^{j_2} = O(1)$ for $-1\le t\le -1+\frac 1 C$. On the other hand, using the given hypothesis we may estimate
\begin{eqnarray*}
{j\choose k}b_{j}- {j-1 \choose k}b_{j-1} 
&=&  {j\choose k}(b_{j}- b_{j-1}) + b_{j-1}{j-1\choose k-1} \\
&=& O((j+1)^{\rho_1+k}) + O((j+1)^{\rho+k-1}).
\end{eqnarray*}
Without loss of generality we may assume $\rho_1>\rho-1$. Since $|t|^k\sim 1$, we obtain
\begin{eqnarray*}
|m^{(k)}_n(t)| &\lesssim&  \sum_j (j+1)^{\rho_1+k} |t|^j  \quad \lesssim \quad (1-|t|+\frac 1 n)^{-(\rho_1+k+1)}\\ &\lesssim& (1+\frac 1 n - |t|)^{\epsilon} \sqrt{Var[r^{(k)}_n(t)]},
\end{eqnarray*}
where $\epsilon = \rho-\rho_1-\frac 1 2>0$. 

Similarly, for $m^*_n$ we may estimate, with the assistance of Lemma~\ref{l.elementary},
\begin{eqnarray*}
(k!)^{-1}{m^*}^{(k)}_n(t)  &=&   b_{n-k}(1+\dots+t^{n-k}) + \\
&& + \ \ \sum_{j=k+1}^{n} \Big({j\choose k}b_{n-j} - {j-1 \choose k}b_{n-j+1}\Big)(t^{j-k}+\dots + t^{n-k}) \\
&\lesssim& \sum_{j} (j+1)^{k} (n-1+j)^{\rho_1} |t|^{j}  + \sum_{j} (j+1)^{k-1} (n+1-j)^{\rho} |t|^{j} \\
&\lesssim& n^{\rho_1} (1 +\frac 1 n - |t|)^{-k-1}  + n^{\rho}(1+\frac 1 n - |t|)^{-k}\\ 
&\lesssim&  (1-|t|+\frac 1 n)^{\epsilon} \sqrt{Var[{r^*}^{(k)}_n(t)]}.
\end{eqnarray*}
Thus Theorem~\ref{t.hyperbolic-gen} follows from Corollary~\ref{c.mixedmean}.

\section{Correlation functions: background and main estimates}\label{s.corr}
In this section, we summarize our main results about correlation functions for $p_n$ and $p^*_n$. These estimates are key ingredients in the proof of Theorem~\ref{t.general} and the proof for these estimates will be presented in subsequent sections.

We first recall some background about correlation functions, following  \cite{tv2015, dnv2017}.  While there is a more general theory of correlation functions for random point processes,  see for instance \cite{hkpv2009},  our discussion will specialize to the context of the roots of random polynomials. Let $Z$ denote the multi-set of the (complex) roots of $p_n$, where  a root of multiplicity $m$  will be identified as $m$ different elements. 

For $k\ge 1$, we say that a Borel measure  $d\sigma$   on $\C^k$  is  the $k$-point  correlation measure  for the (complex) roots of $p_n$  if the following equality holds for any continuous and compactly supported function $\phi:\C^k \to \C$:
$$\E \sum_{\alpha_1,\dots,\alpha_k \in Z} \phi(\alpha_1,\dots,\alpha_k) = \int_{z\in \mathbb C^k} \phi(z)d\sigma(z).$$
Here, the summation on the left hand side (inside the expectation) is over all ordered $k$-tuples of different elements of $Z$.  The existence of such a measure is a simple application of the Riesz representation theorem. In the literature,  it is common (see e.g. \cite{tv2015}) to define the $k$-point correlation function as the density of $d\sigma$ with respect to the Lebesque measure (which exists for instance in  Gaussian settings \cite{hkpv2009} or more generally smooth distributions), here we will work with correlation measures to allow for more generality.

When $p_n$ is a real polynomial (i.e. with  real-valued coefficients), the set of complex zeros for $p_n$ is symmetric with respect to the real line, and there may be a nontrivial probability that $p_n$ has at least one real root.  Thus, for such polynomials  we will define the mixed complex-real correlation measures for the roots   as follows. Let $m \ge 1$ and $k\ge 0$ and let $d\sigma$ be a measure on $\R^m \times (\C\setminus \R)^k$. We say $d\sigma$ is the  $(m,k)$-point   correlation measure for $Z$ if the following two conditions hold:

(i) $d\sigma$ is symmetric under complex conjugations: for any measurable $A\subset \R^m \times (\C\setminus \R)^k$, it holds that $\rho(A)=\rho(A')$ where $A'$ is one of the $k$ sets obtained from $A$ by taking conjugate in  one fixed coordinate;

(ii) for any compactly supported continuous   $\phi:\R^m \times \C^k \to \C$ we have
$$\E \sum_{ \alpha_i \in Z\cap \R} \sum_{\beta_j\in Z\cap \C_+} \phi(\alpha_1,\dots,\alpha_m,\beta_1,\dots,\beta_k) = \int_{(w,z)\in \mathbb R^m \times \mathbb C_+^k} \phi(w,z)d\sigma(w,z). $$
Here, the summations on the left hand side are over ordered tuples of different elements of $Z$. 
If $d\sigma$ has a density   with respect to the Lebesgue measure, such density  is classically called the $(m,k)$-point correlation function \cite{tv2015}, which will then be  invariant under taking complex conjugation of any variable.

We now define the admissible local sets  where   comparison estimates for the correlation measures will be proved. These are sets where   the expected number of complex roots for $p_n$ could be as small  as a bounded constant $O(1)$. For random polynomials with centered-coefficients, the structure of these sets is well-known and has been exploited by previous authors, here we will use the same structure for random polynomials with non-centered coefficients, following \cite{dnv2017}.

Let  $\delta>0$ that may depend on $n$. Define
\begin{eqnarray}\label{e.Idelta}
I(\delta)  =  \begin{cases}\{z\in \C: 1-2\delta \le |z| \le 1-\delta\}, & \delta \ge \frac 1 {10n};\\
\{z\in \C: 1-\frac 1{2n}\le |z|\le 1+\frac 1 {2n}\}, & \delta <\frac 1 {10n};
\end{cases}
\end{eqnarray}
Define
$I_{\R}(\delta) = I(\delta) \cap \R$ and define $I_{\C_+}(\delta)=I(\delta)\cap \C_+$.

Let $p^*_n(z):=z^n p_n(1/z)$ be the reciprocal polynomial of $p_n$.

Below, we say that two (possibly complex valued) random variables $\xi_j$ and $\widetilde \xi_j$ have matching moments to up to second order if
\begin{eqnarray}\label{e.matchingmoment}
\E Re(\xi_j)^\alpha Im(\xi_j)^\beta = \E Re(\widetilde \xi_j)^\alpha Im(\widetilde \xi_j)^\beta
\end{eqnarray}
for any $0\le\alpha,\beta\le 2$ such that $\alpha+\beta \le 2$. Note that if one of $\xi_j$, $\widetilde \xi_j$ is real valued then this matching condition will force the other to be real-valued. The Gaussian analogue of $p_n(z)=\sum_j (b_j+c_j\xi_j)z^j$ if $G_j$ is defined to be  $p_{n,G}(z)=\sum_{j} (b_j+c_j G_j)z^j$ where $G_0,\dots, G_n$ are independent Gaussian and $G_j$ and $\xi_j$ have matching moments  up to the second order.

Our main result about the mixed complex-real $(m,k)$-point correlation functions for the roots of $p_n$ is stated below, here $m\ge 1$ and $k\ge 0$. In Theorem~\ref{t.corr}, we consider a real random polynomials whose coefficients satisfy Condition~\ref{cond.polyvar}, and we let $d\sigma$ and $d\sigma^*$ denote the $(m,k)$-point  correlation measures  for the roots of $p_n$   and $p^*_n$. The Gaussian analogues of these two correlation measures will be denoted by $d\sigma_G$ and $d\sigma^*_G$.  

In the following, it is understood that all implicit constants  may depend on  the implicit constants in Condition~\ref{cond.polyvar}.  


\begin{theorem}\label{t.corr} 

Given $0<c<\widetilde c <1$,  we could find   $C_1,\alpha_1>0$ such that the following holds for any $\frac 1 n \lesssim \delta \le \frac 1 {C_1}$ and any $(x,z)=(z_1,\dots,z_m,z_{m+1},\dots,z_{m+k})   \in I_{\R}(\delta)^{m} \times I_{\C_+}(\delta)^k$:\\

Let   $\phi_\delta$ be supported on $B_{\R}(0,c\delta)^m \times B_{\C}(0,c\delta)^k$ such that as a function on $\R^{m+2k}$ it is in $C^{3k+2}$ and furthermore $\sup|\partial^\alpha \phi_\delta| \le \delta^{-|\alpha|}$ up to order $|\alpha|\le 3k+2$.  \\

Let $J \subset I_{\R}(\delta)+(-\widetilde c \delta, \widetilde c\delta)$ be such that for any $1\le j\le m+k$ the following holds\footnote{Note that the interval $J=I_{\R}(\delta)+(-\widetilde c \delta, \widetilde c\delta)$ has this property, although in the applications we may work with much thinner intervals (which is allowed if $\widetilde c$ is small).}: 
\begin{itemize}
\item if $sign(Re(z_j))\ge 0$ and $|Im(z_j)|\le \widetilde c \delta$ then $(|z_j|-\widetilde c \delta, |z_j|+\widetilde c\delta)\subset J$.  
\item if $sign(Re(z_j))< 0$ and  $|Im(z_j)|\le \widetilde c \delta$ then $(-|z_j|-\widetilde c \delta, -|z_j|+\widetilde c\delta)\subset J$.  \\
\end{itemize} 


(i) Assume that   $|m''_n|\lesssim \sqrt{Var[r''_n]}$ uniformly on $J$, or $|m_n|>C_1|\log(1+\frac 1 n - |t|)|^{1/2}\sqrt{Var[r_n]}$ for all $t\in J$.  Then
\begin{eqnarray*}
&& \int_{\R^m\times \C^k_+} \phi_\delta(y-x, w-z) [d\sigma (y,w)-d\sigma_G(y,w)] = O(\delta^{\alpha_1}) .\\
\end{eqnarray*}

(ii) Assume that   $|{m^*}''_n|\lesssim \sqrt{Var[{r_n^*}'']}$ uniformly on $J$, or $|m^*_n|>C_1|\log(1+\frac 1 n - |t|)|^{1/2}\sqrt{Var[r^*_n]}$ for all $t\in J$.  Then
\begin{eqnarray*}
&&\int_{\R^m\times \C^k_+} \phi_\delta(y-x, w-z) [d\sigma^* (y,w)-d\sigma^*_G(y,w)] = O(\delta^{\alpha_1}).
\end{eqnarray*}
\end{theorem}

Our proof   will use the following   result for  the $k$-point complex correlation functions, where $k\ge 1$. In Theorem~\ref{t.complex-corr}, we consider a (possibly complex valued) random polynomial $p_n$ whose coefficients satisfy Condition~\ref{cond.polyvar}. Below we let $d\sigma$ and $d\sigma^*$ denote the $k$-point correlation measures for the zeros of $p_n$ and $p^*_n$, and let  $d\sigma_G$ and $d\sigma^*_G$ be their Gaussian analogues.

\begin{theorem}\label{t.complex-corr} Given any  $0<c<1$, we could find constants $C_1,\alpha_1>0$ such that the following holds for any $\frac 1 n \lesssim \delta \le \frac 1 {C_1}$ and any $z\in I(\delta)^k$:\\

Let $\phi_\delta$  be supported on $ B_{\C}(0, c\delta)^k$ such that as a function on $\R^{2k}$ it is $C^{3k+2}$ and furthermore $\sup|\partial^\alpha \phi_\delta| \le \delta^{-|\alpha|}$ up to order $|\alpha|\le 3k+2$.

Then
\begin{eqnarray*}
&&\int_{\C^k} \phi_\delta (w-z) [d\sigma (w)-d\sigma_G(w)]  \quad =\quad O(\delta^{\alpha_1}),\\
&& \int_{\C^k} \phi_\delta(w-z) [d\sigma^* (y,w)-d\sigma^*_G(y,w)] \quad =\quad  O(\delta^{\alpha_1}).
\end{eqnarray*}

\end{theorem}

Our Theorem~\ref{t.complex-corr} slightly generalizes   \cite[Theorem 2.3]{dnv2017}. Here we point out an example outside the scope of \cite{dnv2017}. Recall that in \cite[Theorem 2.3]{dnv2017} it is assumed that $p_n(z)=c_0\xi_0+c_1\xi_1 z + \dots + c_n \xi_n z^n$ where $\xi_j$ are independent with unit variance (but could have nonzero means).  In our setting, with $p_n(z)=a_0+a_1z +\dots + a_n z^n$, if  $a_j$ is a nonzero constant with probability $1$ (which is allowed to happen for $j=O(1)$ or $j\ge n-O(1)$ according to Condition~\ref{cond.polyvar}) then it is not possible to write $a_j = c_j \xi_j$ where $\xi_j$ of variance $1$.

We will prove Theorem~\ref{t.complex-corr} using an adaptation of the proof of \cite[Theorem 2.3]{dnv2017}.  We   take this as an opportunity to provide a more streamlined presentation of the argument in \cite{dnv2017}, in particular in the proof we will   prove  new estimates involving log integrability of random polynomials and  bounds on the local number of roots, which could be of independent interests.

\section{Local anti-concentration inequalities}
In this section we will prove several anti-concentration inequalities for  random polynomials whose coefficients  satisfy Condition~\ref{cond.polyvar}. We will use these estimates later in the proof of Theorem~\ref{t.complex-corr}. Below, let $q_n= (n+1)^{-\rho} p^*_n$ be the normalized  reciprocal polynomial for $p_n$. Recall that 
$$I(\delta) = \begin{cases} \{z\in \C: 1-2\delta \le |z|\le 1-\delta\}, &\text{if $\delta \ge \frac 1 {10n}$};\\
\{z\in \C: 1-\frac 1{2n}\le |z|\le 1+ \frac1 {2n}\}, & \text{if $\delta<\frac 1 {10n}$}.
\end{cases}$$

Our first set of estimates is contained the following theorem:
\begin{theorem}\label{t.smallball1} Let $0\le c <1$. Then there are constants $C_1, \alpha_1>0$  such that the following holds for any $\frac 1 n \lesssim \delta \le \frac 1{C_1}$ and  any $|z|\in I(\delta)+(- c\delta, c\delta)$ and any $t>0$:
\begin{eqnarray}
\label{e.smallball1}       \sup_{u} \P(|p_n(z)-u| \le t) &\lesssim&  (t\delta^\rho)^{\alpha_1} +   e^{-\alpha_1 n\delta}, \\
\label{e.smallball1*}      \sup_{u} \P(|q_n(z)-u| \le t) &\lesssim&    t^{\alpha_1} +   e^{-\alpha_1 n\delta}.
\end{eqnarray}
\end{theorem}

Now, if $\delta \approx1/n$ then Theorem~\ref{t.smallball1} does not give us much information:  the right hand sides of \eqref{e.smallball1} and \eqref{e.smallball1*} are now comparable to $1$, therefore these estimates  hold automatically. In this range  of $\delta$,  the following set of estimates is more useful. Below, let $\log_+(x)=\max(0,\log x)$.

\begin{theorem}\label{t.smallball2} Let $0\le c <1$. Then there is a constant $C_1>0$  such that the following holds for any $\frac 1 n \lesssim \delta \le \frac 1 {C_1}$ and any $|z|\in I(\delta) +(-c\delta, c\delta)$ and any $t>0$:
\begin{eqnarray}
\label{e.smallball2}      \sup_u \P(|p_n(z)-u| \le t) &\lesssim&  n^{-1/2}+ \delta^{1/2}\log_+^{-1/2}(\frac 1{t\delta^{\rho}}) , \\
\label{e.smallball2*}   \sup_u \P(|q_n(z)-u| \le t) &\lesssim&   n^{-1/2}+ \delta^{1/2}\log_+^{-1/2}(\frac 1{t}).
\end{eqnarray}
\end{theorem}

As a  corollary of  Theorem~\ref{t.smallball1} and Theorem~\ref{t.smallball2}, we obtain

\begin{corollary} \label{c.logpn}
Let $0\le c <1$. Then there is a constant $C_1>0$   such that the following holds for any $\frac 1 n \lesssim \delta  \le \frac 1 {C_1}$ and any $|z|\in I(\delta)+(-c\delta, c\delta)$: for any  $0<\alpha_2<\frac 1 2$ there is a constant $C_2$ such that
\begin{eqnarray*}
\P(\log |p_n(z)| \le -C_2 |\log \delta|)  &\lesssim& \delta^{\alpha_2}.\\
\P(\log |q_n(z)| \le -C_2 |\log \delta|) &\lesssim& \delta^{\alpha_2}.
\end{eqnarray*}
\end{corollary}
\proof[Proof of Corollary~\ref{c.logpn}] Below we only prove the claimed estimate for $\log|p_n|$, and  the  same argument specialized to the case $\rho=0$ can be applied to $\log|q_n|$. Using Theorem~\ref{t.smallball1} and Theorem~\ref{t.smallball2},  for  any $\lambda > 0$ we have
$$\P(\log |p_n(z)| \le  (\rho-\lambda)|\log \delta|)  \ \ \lesssim \  \ \min(\delta^{\alpha_1\lambda} +   e^{-\alpha_1 n\delta}, 
n^{-1/2}+\lambda^{-1/2} (\frac{\delta}{|\log \delta|})^{1/2}).$$
Thus, for any $\delta \in [\frac{\alpha_2}{\alpha_1}\frac {\log n}n, \frac 1 {C_1}]$ we have
$$P(\log |p_n(z)| \le -(\frac {\alpha_2}{\alpha_1}-\rho)|\log \delta|) \quad \lesssim \quad \delta^{\alpha_2} + e^{-\alpha_2 \log n} \quad \lesssim  \quad \delta^{\alpha_2}.$$
On the other hand, for any $\frac 1 n \lesssim \delta \le \frac {\alpha_2}{\alpha_1}\frac{\log n}n$ we have
$$P(\log |p_n(z)| \le -(\frac {\alpha_2}{\alpha_1}  - \rho)|\log \delta|) \quad \lesssim \quad n^{-1/2}+\delta^{1/2}|\log \delta|^{-1/2} \quad \lesssim \quad \delta^{\alpha_2}.$$
\endproof

\subsection{Proof of Theorem~\ref{t.smallball1}} Recall that $p_n(z) = \sum_j (b_j+c_j \xi_j) z^j$.  Using Condition~\ref{cond.polyvar}, we may find  $j_0 \ge 0$ and $M_0>0$ such that 
\begin{eqnarray}
\label{e.condcj}
|c_j| \le M_0 (1+j)^\rho
\end{eqnarray}
for all $j$, while  $|c_j|\ge M_0^{-1}(1+j)^\rho$ for $j_0\le j\le n-j_0$.

We first prove \eqref{e.smallball1}. Since the left hand side of \eqref{e.smallball1} is $O(1)$,   we may assume without loss of generality that  $\delta>\frac B{n}$ for a large absolute constant $B$. In particular, we will have $1-(2+c)\delta \le |z| \le 1-(1-c)\delta$, thus $ |z|^N \le (1-(1-c)\delta)^{N}$.

Now, there is a  constant $c'>0$ depending only on $c$ such that $(1-(1-c)\delta)^{1/\delta}<1-c'$ for all $\delta>0$. Therefore, we may choose $j_0\le N \approx 1/\delta$ such that $|z|^N$ is very small. In particular, we may choose such $N$ so that  $|z|^N < 2^{-(\rho+2)}M_0^{-2}$. Now, observe that, thanks to \eqref{e.condcj},
$$|c_{kN}/c_{(k+1)N}| \ge 2^{-(\rho+1)}M_0^{-2}$$ 
for any $1\le k \le (n-j_0)/N$. Therefore, 
\begin{eqnarray}
\label{e.clacunary}
|c_{N}z^N| \ge 2 |c_{2N}z^{2N}|\ge \dots \ge 2^{\ell-1} |c_{\ell N}z^{\ell N}|
\end{eqnarray}
for any $1\le \ell \le [\frac{n-j_0}N] \approx n\delta$.

We now recall  the following  anti-concentration bound:

\begin{claim}
Let $\epsilon_0,C_0 >0$. Then there are  constants $\alpha_2, C_2>0$ such that the following holds for any  $\ell\ge 1$: If $\xi_1,\dots, \xi_\ell$ are independent  with zero mean and unit variance satisfying $\E |\xi_j|^{2+\epsilon_0}<C_0$, then for any lacunary sequence $|d_1|\ge 2|d_2|\ge \dots \ge 2^{\ell-1}|d_\ell|$  we have:
$$\sup_u \P(|\sum_{j=1}^m d_j \xi_j -u| \le  |d_\ell|) \le C_2 e^{-\alpha_2\ell}.$$ 
\end{claim}

For a proof of this now-standard bound, see e.g. \cite[Lemma 9.2]{tv2015} or \cite[Lemma 4.2]{dnv2017}. We apply the above anti-concentration bound to $d_j=c_{jN}z^{jN}$ and to the random variables $\xi_{N},\dots, \xi_{(\ell-1)N}$. By absorbing the remaining   terms in $p_n(z)$  into the concentration point $u$, it follows that
\begin{eqnarray}\label{e.presmallball1}
\sup_u \P(|p_n(z)-u| \le |c_{\ell N}z^{\ell N}|) = O(e^{-\alpha_2\ell}),
\end{eqnarray}
for any $1\le \ell \le \ell_N:= [(n-j_0)/N]$. To obtain the desired estimate \eqref{e.smallball1} from this inequality, we will choose $\ell$ to depend on $t$, and this choice is explained below.

First, note that $|z|^{1/\delta} \ge (1-(2+c)\delta)^{1/\delta}$, which is uniformly bounded away from $0$ and since $N\approx 1/\delta$, we may find  a constant  $\alpha_3>0$ such that  $ |z^N| \ge e^{-\alpha_2/2}$. 
It follows that
$$|c_{\ell N} z^{\ell N}| \gtrsim  (\ell N)^\rho e^{-\alpha_3\ell/2} \gtrsim  N^\rho e^{-\alpha_3\ell} \gtrsim \delta^{-\rho}e^{-\alpha_3\ell}$$
For convenience, let $C_3>0$ be such that $|c_{\ell N} z^{\ell N}| \ge \frac 1 {C_3} \delta^{-\rho} e^{-\alpha_3\ell}$.
We then let  $\ell$  to be the integer such that
$$\frac 1 {C_3}  e^{-\alpha_3(\ell+1)} \le  t\delta^{\rho}<\frac 1 {C_3}  e^{-\alpha_3\ell}.$$

Now, since the left hand side of \eqref{e.smallball1} is $O(1)$ we may assume without loss of generality that $\ell \ge 1$. To check that this $\ell$ will lead us to \eqref{e.smallball1}, we divide the consideration into two cases:

\emph{\underline{Case 1:}  $1\le \ell \le \ell_N$.} 

It follows from the above constraint on $\ell$ that  $e^{-\ell} = O((t\delta^\rho)^{1/\alpha_3})$. In this range of $\ell$ we may use  \eqref{e.presmallball1}, and obtain
\begin{eqnarray*}
\sup_u \P(|p_n(z)-u|\le t)  &\le& \P(|p_n(z)| \le \frac 1 {C_3}\delta^{-\rho} e^{-\alpha_3\ell})\\
&\le&  \sup_u \P(|p_n(z)-u|\le |c_{\ell N} z^{\ell N}|) \\
&\lesssim&  e^{-\alpha_2\ell}  = O((t\delta^\rho)^{\alpha_2/\alpha_3}).
\end{eqnarray*}
Thus  by ensuring $\alpha_1\le \alpha_2/\alpha_3$ we obtain \eqref{e.smallball1}.

\emph{\underline{Case 2:}  $\ell > \ell_N$.}  

Here \eqref{e.presmallball1} is not available, however   we observe that  the LHS of \eqref{e.smallball1}  is  nondecreasing with respect to $t$. Therefore, using the case $\ell=\ell_N$ of Case 1, we obtain
\begin{eqnarray*}
\sup_u \P(|p_n(z)-u|\le t)  &\le&   \sup_u \P(|p_n(z)-u|\le |c_{\ell_N} z^{\ell_N N}|)  \ \lesssim \  e^{-\alpha_2\ell_N}.
\end{eqnarray*}
Since $\ell_N  \approx n\delta$, the last estimate can be bounded above by $O(e^{-\alpha_1 n\delta})$ for some $\alpha_1>0$.
This completes the proof of  \eqref{e.smallball1}.

We now discuss the proof of \eqref{e.smallball1*}, which will follow the same argument. For convenience of notation, we let $q_n(x)=(e_0+d_0 \widetilde \xi_0) + (e_1+d_1 \widetilde \xi_1) x+\dots + (e_n+d_n \widetilde \xi_n) x^n$, where $e_j = b_{n-j} (n+1)^{-\rho}$, $d_j=c_{n-j} (n+1)^{-\rho}$ and $\widetilde \xi_j = \xi_{n-j}$. It is clear that $e_j \lesssim 1$ and $d_j \approx 1$ for $j_0\le j\le n/2$, therefore we may apply  the special case $\rho=0$ of \eqref{e.smallball1}   to the random polynomial $d_0\widetilde \xi_0 + \dots + d_{[n/2]} \widetilde \xi_{[n/2]} x^{[n/2]}$. The desired estimate for $q_n$ then follows by absorbing the other terms into the concentration point $u$.

\subsection{Proof of Theorem~\ref{t.smallball2}}

Below we only prove \eqref{e.smallball2}, and \eqref{e.smallball2*} can be obtained from \eqref{e.smallball2} by arguing as in   the proof of Theorem~\ref{t.smallball1} in the last section. 

The proof   uses the following generalization of a lemma of Erd\"os (for a proof see \cite[Lemma 4.1]{dnv2017}):

\begin{claim}Let $\epsilon_0,C_0>0$. Then there is a constant $C>0$ such that the following holds for any $m\ge1$: If  $\xi_1,\dots, \xi_m$ are independent   and $\sup_j \E |\xi_j|^{2+\epsilon_0}<C_0$ then for any $d_1,\dots, d_m\in \C$ we have
$$\sup_u \P(|d_1 \xi_1+\dots+d_m\xi_m-u| \le \min|d_j|/C)\le C/\sqrt m.$$
\end{claim}
Let $n-j_0\ge m\ge 2j_0$, where $j_0=O(1)$ is such that  $|c_j|$ is comparable to $(1+j)^\rho$ for $j_0\le j\le n-j_0$ (thanks to Condition~\ref{cond.polyvar}).  Applying the above estimate to $d_j=  c_j z^j$ for $m/2\le j\le m$, it follows that
$$\sup_u \P(|p_n(z)-u| \le \min_{m/2\le j \le m} |c_jz^j|/C) = O(1/\sqrt m)$$
Now, we may choose  $C\ge 1$ be sufficiently large such that  $\delta \ge 1/(C n)$. For any $z\in I(\delta)+(-c\delta, c\delta)$, it holds that $|z|\ge 1- 2C\delta$, therefore
$$\min_{m/2\le j\le m} |c_jz^j| \gtrsim  m^\rho  (1-2C\delta)^m \gtrsim m^\rho e^{-2Cm\delta} \gtrsim \delta^{-\rho}  e^{-3Cm\delta}.$$ 
Collecting estimates, for   $C>0$ large enough we will have
\begin{eqnarray}
\label{e.presmallball2}
\sup_{u} \P(|p_n(z)-u| \le \frac 1 {C} \delta^{-\rho} e^{-Cm\delta}) \quad =\quad O(m^{-1/2}),
\end{eqnarray}
for any integer $m\in [2j_0, n-j_0]$. To obtain the desired estimate \eqref{e.smallball2} from this inequality, we will choose $m$ suitably depending on $t>0$. We will choose $m$ to be the integer such that 
$$\frac 1 {C}  e^{-C(m+1)\delta} < t\delta^\rho \le \frac 1{C} e^{-Cm\delta}.$$
Now, since the LHS of \eqref{e.smallball2} is $O(1)$, we may assume without loss of generality that $m\ge 2j_0$.  To show that this choice would give us \eqref{e.smallball2}, we divide the consideration into  two cases:

\underline{\emph{Case 1}}: $2j_0\le m\le n-j_0$. For such $m$ we may use \eqref{e.presmallball2}. We note that, as a consequence of the above constraint on $m$, we will have $m\delta\gtrsim   \log_+(\frac 1 {t\delta^\rho}) $. Consequently,
$$\sup_u \P(|p_n(t) - u | \le t) \quad \lesssim \quad m^{-1/2}  \quad \lesssim \quad \delta^{1/2}\log_+^{-1/2}(\frac 1{t\delta^{\rho}}).$$

\underline{\emph{Case 2}}:  $m \ge n-j_0+1$. Here we will use monotonicity of the left hand side of \eqref{e.smallball2} (as a function of $t$). Since we now have $t < \frac 1 {C}\delta^{-\rho} e^{-C(n-j_0)\delta}$, it follows that
\begin{eqnarray*}
\sup_u \P(|p_n(z)-u| \le t) &\le& \sup_u \P(|p_n(z)-u| \le \frac 1 {C} \delta^{-\rho} e^{-C(n-j_0)\delta})  \ \lesssim \ n^{-1/2}.
\end{eqnarray*}
This completes our proof of Theorem~\ref{t.smallball2}.

\section{Logarithmic integrability of random polynomials}\label{s.loginteg}

This section is devoted to establishing several estimates about the  integrability of  $\log|p_n|$ and $\log |p^*_n|$, which will be used to prove bounds for the number of local real roots of $p_n$ in subsequent sections.  Throughout this section, we'll assume that the coefficients of $p_n$ satisfy Condition~\ref{cond.polyvar}. For convenience, let $q_n:=(n+1)^{-\rho} p^*_n$. 

\subsection{Logarithmic integrability on the unit disk}

 We start with an estimate about integrability  on the unit disk $B(0,1)=\{|z|\le 1\}$. We view this  as a global estimate. 

\begin{theorem}\label{t.nns-weak} There are absolute constants $C,c>0$ and an event $F$ of exponentially decaying probability $\P(F)=O(e^{-c n})$   such that   the following holds:  
\begin{eqnarray}
\label{e.nns-weak} \E [1_{F^c} \int_{B(0,1)} |\log|p_n(w)||^q  dw] &\le&  (C q)^{C q}   (\log (n+2))^{Cq} 
\end{eqnarray}
for all $q\ge 1$, and the analogous estimate  also holds for $q_n$. 
\end{theorem}

We note that the exclusion of  an exceptional set   of exponentially decaying probability is  important. To see this, suppose that  $b_j=0$ for all $j$, then $p_n(x)\equiv 0$ on  the event $F=\{\xi_j=0 \ \ \forall j\}$, which has an exponentially decaying probability $\P(F)=O(p^{n})$ if  for some fixed $p\in (0,1)$ we have $\P(\xi_j=0)\ge p$ for all $j$. Such event must be excluded to ensure any integrability for $|\log |p_n||$ on $B(0,1)$.

Without loss of generality we   may assume that $n\ge 3$ in the proof.  Given such a condition,   the right hand side of  \eqref{e.nns-weak} is a strictly increasing function of the implicit constant $C$, which will be convenient in the proof.

To start, we note that the estimate \eqref{e.nns-weak}  follows from a slightly weaker estimate:
\begin{proposition}\label{t.nns-weak1} There is an event $F$ of exponentially decaying probability $\P(F)=O(e^{-c n})$ (for some fixed $c>0$) such that   the following holds: for any $\epsilon>0$, there is a constant $C=C(\epsilon)$  such that
\begin{eqnarray}
\label{e.nns-weak1} \E [1_{F^c} \int_{B(0,1)} |\log|p_n(w)||^q  dw] &\le&  (C q)^{C q} n^{C} (\log (n+2))^{(1+\epsilon) q} 
\end{eqnarray}
for all $q\ge 1$, and the analogous estimate  also holds for $q_n$. 
\end{proposition}
Indeed,  the key observation here is that the  the implicit constant $C$ does not depend on $q$.  If \eqref{e.nns-weak1} holds, using Holder's inequality we have, for any $p\ge 1$:
\begin{eqnarray*}
 \E [1_{F^c} \int_{B(0,1)} |\log|p_n(w)||^q  dw] &\le& \Big(\E \Big[1_{F^c} \Big(\int_{B(0,1)} |\log|p_n(w)||^q  dw\Big)^p\Big]\Big)^{1/p}\\
&\le& \pi^{1-\frac 1 p} \Big(\E [1_{F^c}  \int_{B(0,1)} |\log|p_n(w)||^{pq}  dw]\Big)^{1/p}\\
&\le& \pi^{1-\frac 1 p} \Big((Cpq)^{Cpq}n^C (\log n)^{(1+\epsilon)pq}\Big)^{1/p}\\
&=& \pi^{1-\frac 1 p}  (Cpq)^{Cq}n^{C/p} (\log n)^{(1+\epsilon)q} 
\end{eqnarray*}
The desired conclusion \eqref{e.nns-weak} then follows by choosing $p=\log n$.

The main ingredient in the proof of Proposition~\ref{t.nns-weak1} is a result of Nazarov-Nishry-Sodin \cite[Corollary 1.2]{nns2014} for random Fourier series, summarized below:
\begin{proposition}\label{p.nns}\cite{nns2014} There is an absolute constant $C>0$ such that the following holds: Let $r_{\epsilon}(z)=\sum_j \epsilon_j d_j z^j$ where $d_j$ are deterministic with $\sum_j |d_j|^2=1$ and $\epsilon_j$ are independent Rademacher random variables. Then for any $p>0$
$$\E [\int_{0}^{2\pi} \int_0^1  |\log|r_{\epsilon}(re^{i\theta})||^p drd\theta] \le (C p)^{6p}.$$
\end{proposition}

Our proof will actually use the following simple  extension of Proposition~\ref{p.nns}.

\begin{lemma}\label{l.nns-deviation}There is an absolute constant $C>0$ such that the following holds for any $m:B(0,1)\to \mathbb C$   measurable with $M:=\int_0^{2\pi} \int_0^1 |m(re^{i\theta})|^2 drd\theta <\infty$: Let $r_{\epsilon}(z)=\sum_j \epsilon_j d_j z^j$ where $d_j$ are deterministic with $\sum_j |d_j|^2=1$ and $\epsilon_j$ are independent Rademacher random variables.  Then for any $p>0$   we have 
$$\E \int_0^{2\pi}\int_0^1 |\log |m(re^{i\theta})+r_{\epsilon}(re^{i\theta})||^p drd\theta \lesssim (C p)^{7p} (M+1).$$
\end{lemma}
In Lemma~\ref{l.nns-deviation},  we could in fact replace  the constant $7$  by any constant bigger than $6$ (for our applications any absolute constant would suffice).  

\subsubsection{Proof of Lemma~\ref{l.nns-deviation}}

To prove Lemma~\ref{l.nns-deviation}, we will use the following  crude estimate. For convenience of notation, let $f(z)=m(z)+r_\epsilon(z)$ and let $|.|$ denote the Lebesgue measure of measurable subsets of $[0,1]\times [0,2\pi]$.

\begin{claim}\label{cl.nns-deviation2} There is an absolute constant $C>0$ such that for any $p>0$ and $\lambda\ge 0$ we have 
\begin{eqnarray*}
\E |\{(r,\theta): \log |f(re^{i\theta})| > \lambda\}|   &\lesssim& (1+\lambda)^{-p} (Cp)^{p} (M+1).
\end{eqnarray*}
\end{claim}
\begin{align*} \text{Indeed,} &&& \E |\{(r,\theta): \log |m(re^{i\theta})+r_{\epsilon}(re^{i\theta})| > \lambda\}| \\
&\lesssim&&  e^{-2\lambda}\Big( \E \int\int |m(re^{i\theta})|^2 drd\theta +  \E \int\int |r_{\epsilon}(re^{i\theta})|^2 drd\theta\Big)\\
&\lesssim&&  e^{-2\lambda} (M + \int\int \sum_j |d_j(re^{i\theta})^j|^2 drd\theta)\\
&\lesssim&&  e^{-2\lambda} (M+1).
\end{align*}
Now, let $h\ge 1$ be integer such that   $h-1<p\le h$, we then have
\begin{eqnarray*}
e^{2\lambda} &\ge& (1+2\lambda)^h/h! > (1+\lambda)^h h^{-(h-1)} \\
&\ge& (1+\lambda)^p (p+1)^{-p} \gtrsim (1+\lambda)^{p} p^{-p}.
\end{eqnarray*}
This competes the proof of Claim~\ref{cl.nns-deviation2}.

In the proof of Lemma~\ref{l.nns-deviation}, we will  use another estimate, which in turn is a consequence of  Proposition~\ref{p.nns}.

\begin{claim}\label{cl.nns-deviation1}There is an absolute constant $C$ such that for any $p>0$ and $\lambda\ge 0$   we have 
\begin{eqnarray*}
\E |\{(r,\theta): \log |f(re^{i\theta})| < - \lambda\}|  &\lesssim& (1+\lambda)^{-p} (C p)^{6p} .
\end{eqnarray*}
\end{claim}

Since the left hand side of the above estimate is always bounded above by $2\pi$ and since $p^{p} \ge e^{-1/e}$ for any $p>0$,  we may assume   $\lambda>1$ without any loss of generality. For such $\lambda$, it suffices to show that
$$\E |\{(r,\theta): \log |f(re^{i\theta})| < - \lambda\}| \lesssim (\lambda-\frac 1 2\ln 2)^{-p} (C p)^{6p} .$$

 Let $\epsilon'_j$ be iid copies of $\epsilon_j$, such that $\epsilon'_0,\dots,\epsilon'_n,\epsilon_0,\dots,\epsilon_n$ are independent Rademacher random variables. Let $\eta_j = (\epsilon_j - \epsilon_j')/\sqrt 2$, which  are also independent Rademacher random variables. We have
\begin{eqnarray}
\nonumber &&\Big( \P(\log |f(re^{i\theta})| < - \lambda)\Big)^2 \\
\nonumber &=&  \P(|m(re^{i\theta})+r_{\epsilon}(re^{i\theta})| < e^{- \lambda}, \ |m(re^{i\theta})+r_{\epsilon'}(re^{i\theta})| < e^{- \lambda}) \\
\nonumber &\le&  \P(|r_{\epsilon}(re^{i\theta})-r_{\epsilon'}(re^{i\theta})| < 2e^{- \lambda}) \\
\label{e.radlog} &\le& (\lambda - \frac 1 2\ln 2)^{-2p} \E |\log |r_{\eta}(re^{i\theta})||^{2p}.
\end{eqnarray}
\begin{align*}
\text{Thus,} &&&\E|\{(r,\theta): |\log |f(re^{i\theta})| < - \lambda\}|  \\
&=&& \int_{0}^{2\pi} \int_0^1  \P(\log |f(re^{i\theta})| < - \lambda) drd\theta  \quad \text{ (by Fubini's theorem)}\\
&\lesssim&& \Big(\int_{0}^{2\pi} \int_0^1  \Big(\P(\log |f(re^{i\theta})| < - \lambda)\Big)^2 drd\theta\Big)^{1/2} \quad \text{ (by H\"older's inequality)}\\
 &\lesssim&& (\lambda - \frac 1 2\ln 2)^{-p}  \Big(\int_0^{2\pi} \int_0^1 \E |\log |r_{\eta}(re^{i\theta})||^{2p} drd\theta\Big)^{1/2} \quad \text{(by \eqref{e.radlog})}\\
 &\lesssim&& (\lambda - \frac 1 2\ln 2)^{-p} (Cp)^{6p} \quad \text{(using Proposition~\ref{p.nns} with $2p$  and choosing a large $C$).}
\end{align*}
This completes the proof of Claim~\ref{cl.nns-deviation1}.

We are now ready to start the proof of Lemma~\ref{l.nns-deviation}. We combine Claim~\ref{cl.nns-deviation1} and Claim~\ref{cl.nns-deviation2} and estimate
\begin{eqnarray*}
&&\E \int_0^{2\pi}\int_0^1 |\log |f(re^{i\theta})||^p drd\theta \\
 &=& p\int_0^\infty  \lambda^{p-1} \E |\{(r,\theta): |\log |f(re^{i\theta})|| > \lambda\}| d\lambda \\
 &\lesssim& (M+1)p\int_0^\infty \lambda^{p-1}(1+\lambda)^{-7p/6} (Cp)^{7p} d\lambda \\
 &\lesssim&(Cp)^{7p} (M+1) \int_0^\infty p(1+\lambda)^{-(1+p/6)}d\lambda \\
 &\lesssim& (Cp)^{7p} (M+1).
 \end{eqnarray*}
This completes the proof of Lemma~\ref{l.nns-deviation}.

\subsubsection{Proof of Proposition~\ref{t.nns-weak1}}
We now start the proof of \eqref{e.nns-weak1}  for $\log |p_n|$. For convenience of notation, we denote $p_{n,\xi}(w)=\sum_{j} (b_j+c_j\xi_j)w^j$ to keep track of the dependence of $p_n$ on the   vector of coefficients  $\xi=(\xi_0,\dots,\xi_n)$. Let
$$F_{\xi}= \{\sigma(\xi) < n^{-1}\},  \text{ where }\sigma(\xi)=(\sum_j |c_j\xi_j|^2)^{1/2}.$$
We first show that $\P(F) = O(e^{-c n})$ for some $c>0$. Since $\xi_j$ are independent and $|c_j|\approx  j^\rho \gtrsim  n^{-1/2}$ for $n-O(1) \ge j\ge O(1)$, it suffices to show that that there are  constants $\delta_0>0$ and $p_0>0$ such that $\P(|\xi_j|<\delta_0) \le 1-p_0$ for all $j$. This was  proved in  Lemma~\ref{l.unif-smallball}.  


We now divide the remaining of the proof  into two cases: the simpler case when $\xi_j$ are symmetric for each $j$, and the general case where no symmetry is assumed. 

\noindent {\bf Case 1: Symmetric coefficients}.

Assume that for each $j$ the distributions of $\xi_j$ and $-\xi_j$ are the same.

Let $\epsilon_0,\dots, \epsilon_n$ be independent Rademacher random variables that are independent from $\xi_0,\dots,\xi_n$, and let $\widetilde \xi_j =\epsilon_j \xi_j$. Thanks to symmetry,  $p_{\xi,n}$ has the same distribution as $p_{\widetilde \xi,n}$. Note that  $\sigma(\xi)=\sigma(\widetilde \xi)$, therefore $F_{\widetilde \xi} = F_\xi$ and is independent of $\epsilon_j$. 
Thus it suffices to show that, for any $C>0$ large enough,
$$\E _{\xi,\epsilon} [1_{F^c_{\xi}} \int_{B(0,1} |\log |p_{n,\widetilde \xi}(w)||^q dw] \lesssim (Cq)^{C q} n^{C} (\log n)^{q} .$$

Note that on the event $F^c_\xi$ we have $\sigma(\xi) \ge n^{-1}$, which implies $|\log \sigma(\xi)|  < \log (n^{2}\sigma(\xi))$. Conditioning on this event and using Lemma~\ref{l.nns-deviation}, we obtain
\begin{align*}
&&&\E_{\epsilon} [  \int_{B(0,1)} |\log |p_{n,\widetilde \xi}||^q] =   \E_{\epsilon}[\int_{B(0,1)}|\log |m_n(w)+\sum_{j} c_j \xi_j \epsilon_jw^j||^q dw] \\
&\lesssim&&  2^q  \E_{\epsilon}\Big[|\log \sigma(\xi)|^q + \int_{B(0,1)}|\log |\frac{m_n(w)}{\sigma(\xi)}+\sum_{j} \frac{c_j \xi_j}{\sigma(\xi)} \epsilon_jw^j||^q dw\Big]\\
&\lesssim&&   \Big[2^q  \log^q (n^2\sigma(\xi)) +  (C p)^{7p} (M+1) \Big]\\
\text{where} &&& M=\int_0^{2\pi} \int_0^1 |\frac{m_n(re^{i\theta})}{\sigma(\xi)}|^2 drd\theta \quad \lesssim  \quad n^2 \sup_{w\in B(0,1)} |m_n(w)|^2 \lesssim n^{C},
\end{align*}
here $C$   depends  on $\rho$. Thus, it remains to show that
$$\E_{\xi} [\log^q (n^2\sigma(\xi))] \le (Cq)^q n^{C} \log^q (n+1)$$
for some $C>0$ (independent of $q$). This estimate in turn follows from concavity of $\log^q(x)$ on $(e^q,\infty)$ and Jensen's inequality:
\begin{eqnarray*}
\E_{\xi} [\log^q (e^{q}+n^2\sigma(\xi))]  &\le& \log^q(\E[e^{q}+n^2\sigma(\xi)]) \\
&\lesssim& \log^q (e^{q}+n^{C}) \quad \lesssim \quad (Cq)^q (\log n)^q .
\end{eqnarray*}

\noindent {\bf Case 2:  General coefficients.}

We now drop the assumption that the distribution of $\xi_j$'s are  symmetric.
To show \eqref{e.nns-weak}, it suffices to show that, for $C=C(\epsilon)>0$ large enough, 
\begin{eqnarray}
\label{e.nns-levelset}
&& \int_{B(0,1)} \P(F_{\xi}^c\cap \{\ell \le |\log |p_{n,\xi}(w)|| \le \ell + 1\}) dw  \\
\nonumber &  \lesssim& (1+\ell)^{-q} (Cq)^{Cq} n^{C} (\log n)^{(1+\epsilon)q} 
\end{eqnarray}
 for any $\ell \ge 0$ and any $q\ge 1$. Since the left hand side of \eqref{e.nns-levelset} is $O(1)$, this estimate holds trivially for $\ell=O(1)$. Thus, we will assume below that $\ell\ge 1$, in particular we may replace $(1+\ell)^{-q}$ by $\ell^{-q}$  on the right hand side without any loss of generality.

Now, let $c'=c/(2q)$. We divide the proof  of \eqref{e.nns-levelset} into two parts, depending on  whether $\ell\le e^{c' n}$ or  $\ell \ge e^{c' n}$.

\underline{\emph{Smaller $\ell$'s:}} For $\ell\le e^{c' n}$, we have $\ell^{-q} \ge e^{-cn/2}$, thus it suffices to show that 
\begin{eqnarray}
\label{e.smalll}
\int_{B(0,1)} \P(|\log |p_{n,\xi}(w)|| \ge \ell) dw  \quad \lesssim \quad  e^{-c n/2}+ \ell^{-q} (C q)^{C q} n^{C} (\log n)^{q}.
\end{eqnarray}
Now, $\{|\log |p_{n,\xi}(w)|| \ge \ell)\} = \{\log |p_{n,\xi}(w)| \ge \ell)\} \cup \{\log |p_{n,\xi}(w)| \le -\ell)\}$, and
\begin{eqnarray*}
\int_{B(0,1)} \P(\log |p_{n,\xi}(w)| \ge \ell) dw 
 &\lesssim&  e^{-2\ell}  \int_{B(0,1)}   \E |p_{n,\xi}(w)|^2 dw \\
&\lesssim& e^{-2\ell} n^{C} \quad \lesssim  \quad \ell^{-q} (Cq)^q n^{C}.
\end{eqnarray*}
Thus, it remains to show that $\int \P(\log |p_{n,\xi}(w)| \le - \ell)$ is bounded by the right hand side of \eqref{e.smalll}.

Let $\widetilde \xi_j$ be iid copy of $\xi_j$ that are   independent of each other and of other $\xi_j$'s. Let $\eta_j = \frac 1 {\sqrt 2} (\xi_j-\widetilde \xi_j)$, then $\eta_j$ is symmetric with mean zero and variance $1$. We also have $\E |\eta_j|^{2+\epsilon_0}=O(C_0)$ uniform over $j$, thanks to Condition~\ref{cond.polyvar}.  One could easily show that $\P(F_{\eta})= O(e^{-cn})$ (with the same $c$ as in the estimate for $\P(F_\xi)$, although this it not important - we could refine the constant $c$ for $F_\xi$ so that    these two exceptional sets share the same constant from the beginning of the proof). 

Now, using H\"older's inequality, we obtain
\begin{eqnarray*}
&& \int_{B(0,1)} \P(\log |p_{n,\xi}(w)| \le -\ell) dw \\
&\lesssim& (\int_{B(0,1)} \P(\log |p_{n,\xi}(w)|, \log |p_{n,\widetilde \xi}(z)| \le -\ell) dw)^{1/2}\\
 &\lesssim&  (\int_{B(0,1)}  \P(\log |p_{n,\eta}(w)| \le -\ell + \frac 1 2 \ln 2) dw)^{1/2}\\
 &\lesssim&  e^{-cn/2}  + (\int_{B(0,1)}  \P(F_{\eta}^c\cap \{\log |p_{n,\eta}(w)| \le -\ell + \frac 1 2 \ln 2\}) dw)^{1/2}\\
 &\lesssim& e^{-cn/2}  + (\ell-\frac 1 2 \ln 2)^{-q} \Big(\E [1_{F_\eta^c} \int_{B(0,1)} |\log |p_{n,\eta}(w)|^{2q}dw] \Big)^{1/2}.
\end{eqnarray*}
Let $C$ be sufficiently large, then using the known estimates for the symmetric case, which applies to $p_{\eta,n}$ and $2q$, we may generously estimate the last display by
\begin{eqnarray*}
&\lesssim& e^{-cn/2}+ (\ell/2)^{-q}  (2Cq)^{Cq} n^{C/2} (\log n)^{q}.
\end{eqnarray*}
This completes the proof of \eqref{e.nns-levelset} for this range of $\ell$.

\underline{\emph{Larger $j$'s:}} For $\ell \ge e^{c'n}$, we proceed as follows. Let $\epsilon_0,\dots, \epsilon_n$ be independent Rademacher random variables that are independent from $\xi_j$'s.  Let $\widetilde \xi_j = \epsilon_j \xi_j$ and consider the symmetrized variant of $p_{n,\xi}$, namely $$p_{n,\widetilde \xi}(z):= \sum_j    (b_j+ c_j\epsilon_j \xi_j)z^j$$

Using H\"older's inequality, for any $p,q \ge 1$ we have
\begin{eqnarray*}
&& \int_{B(0,1)} \P(F_{\xi}^c\cap \{ \ell \le |\log |p_{n,\xi}(w)|| \le  \ell+1\}) dw\\
&\lesssim& \ell^{- (1+\epsilon)q} \Big(\E_{\xi}    [1_{F_{\xi}^c}  \int_{B(0,1)} |\log |p_{n, \xi} (w)||^{(1+\epsilon)pq} dw]\Big)^{1/p} \\
&\le& \ell^{-(1+\epsilon)q} 2^{(n+1)/p}(\E_{\xi}   \E_{\epsilon} [1_{F_{\xi}^c}  \int_{B(0,1)} |\log |p_{\widetilde \xi,n} (w)||^{(1+\epsilon)pq} dw])^{1/p}
\end{eqnarray*}
Here, in the last estimate we used the fact that $p_{n,\xi}$ is equal to $p_{n,\widetilde \xi}$ with probability $2^{-(n+1)}$. Observe that $F_\xi=F_{\widetilde \xi}$. Thus, using the (known) estimate for the symmetric case, we can further estimate the last display by
\begin{eqnarray*}
&\lesssim& \ell^{- (1+\epsilon)q} 2^{n/p} \Big((Cpq)^{Cpq} n^C (\log n)^{(1+\epsilon)pq}\Big)^{1/p}\\
&=& \ell^{- (1+\epsilon) q}  2^{n/p} (Cpq)^{Cq} n^{C/p} (\log n)^{(1+\epsilon)q}.
\end{eqnarray*}
Since $\ell^{\epsilon q} \ge e^{c'nq\epsilon}=e^{cn\epsilon/2}$, it follows that by taking $p\ge \max(1, (c\epsilon)^{-1} \ln 4 )$ we have $\ell^{-\epsilon q} 2^{n/p}\le 1$ and we obtain the desired estimate.

This completes the proof of the desired estimate \eqref{e.nns-weak1} for $\log |p_n|$ of Proposition~\ref{t.nns-weak1}.

We now discuss the proof for the analogous estimate for $\log |q_n|$.  For convenience of notation, let $p^*_n(x) = \sum_{j} (b^*j  + c^*_j \widetilde \xi_j)x^j$ where  $b^*_j=b_{n-j}$, $c^*_j = c_{n-j}$, and $\widetilde \xi_j = \xi_{n-j}$. In particular, $m^*_n(x)=\sum_j b^*_j x^j$. We similarly let 
$$F^*_\xi=\{(n+1)^{-\rho}\sigma^*(\xi) < n^{-1}\}$$ 
where $\sigma^*(\xi) = (\sum_{j} |c^*_j \widetilde \xi_j|^2)^{1/2}$. Using Condition~\ref{cond.polyvar}, we have $|c^*_j|\approx (n+1)^\rho$ for $O(1)\le j\le n/2$, therefore by the same argument as before we obtain $\P(F^*_\xi) = O(e^{-cn})$ for some $c>0$. Now, the proof of the symmetric case is entirely the same as before once we verify that on $F^*_\xi$ it holds that
$$\int_0^{2\pi} \int_0^1 |\frac{m^*_n(re^{i\theta})}{\sigma^*(\xi)}|^2 drd\theta = O(n^C).$$
 But this is clear using Condition~\ref{cond.polyvar}. Finally, the proof of the general case follows from the symmetric case as long as we could verify that $\int_{B(0,1)} \E |q_n(w)|^2 = O(n^C)$, which again is clear from Condition~\ref{cond.polyvar}.

\subsection{Logarithmic integrability on local sets}

In this section we will prove a probabilistic upper bound regarding the local integrability of $\log |p_n|$ and $\log |q_n|$ where $q_n= (n+1)^{-\rho}p^*_n$. This is an estimate on a ball of radius comparable to the scale $\delta$ with center near $I(\delta)$.  All implicit constants below may depend on the implicit constants in Condition~\ref{cond.polyvar}.

\begin{theorem}\label{t.locallogint}  Let $0\le c,c'<1$ be such that $c+c'<1$ and let $C_1>0$ be big enough depending on $c,c'$.  Then for any $\alpha_0\in (0,1/2)$  and $\frac 1 n \lesssim \delta \le \frac 1 {C_1}$ and  $z\in I(\delta)+(-c\delta, c\delta)$ there is an event $F$ with probability $O(\delta^{\alpha_0})$ such that the following estimate holds uniformly  over $1\le p<\infty$:
\begin{eqnarray*}
1_{F^c} \int_{B(z,c'\delta)} |\log |p_n(w)||^p dw  &\le& (C p)^p \delta^{2}   |\log \delta|^{2p}, 
\end{eqnarray*}
and the analogous estimate also holds if we replace $p_n$ by $q_n = (n+1)^{-\rho} p^*_n$.
\end{theorem}
As a consequence  Theorem~\ref{t.locallogint}, we obtain
\begin{eqnarray*}
\E [1_{F^c}\int_{B(z,c'\delta)} |\log |p_n(w)||^p dw ] &\le& (Cp)^p \delta^{2} |\log \delta|^{2p},
\end{eqnarray*}
(and the analogous estimate for $q_n$), which is reminiscent of Theorem~\ref{t.nns-weak}. 

Using Lemma~\ref{l.elementary}, we have the following probabilistic estimates for $\log |p_n|$:

\begin{lemma}\label{l.unif-upper} Let $0\le c<1$. For  $\frac 1 n \lesssim \delta < \frac 1 5$   it holds for any $\epsilon>0$ and $s\in \R$ that 
\begin{eqnarray*}
\P(\sup_{|w|\in I(\delta)+ (-c\delta,  c\delta)} \log |p_n(w)| > s) &\lesssim_\epsilon& e^{-2s}\delta^{-2(\rho+1+\epsilon)}, \\
\P(\sup_{|w|\in I(\delta)+ (-c\delta,  c\delta)} \log |q_n(w)| > s) &\lesssim_\epsilon& e^{-2s}\delta^{-2(1+\epsilon)}.
\end{eqnarray*}
\end{lemma}

\proof[Proof of Lemma~\ref{l.unif-upper}] Recall that $p_n(w)=a_0+a_1w+\dots + a_n w^n$ and $\E |a_j|^{2} = O((1+j)^{2\rho})$ thanks to Condition~\ref{cond.polyvar}. 
Using Lemma~\ref{l.elementary} and Cauchy-Schwartz, for any $\epsilon>0$ and $w\in I(\delta) + (-c\delta, c\delta)$ we have 
\begin{eqnarray*}  
 |p_n(w)|  &\lesssim& (\sum_{j=0}^{n} (1+j)^{2\rho+1+2\epsilon}|w|^{2j})^{1/2} (\sum_{j=0}^n (1+j)^{-2\rho-1-2\epsilon} |a_j|^2)^{1/2}\\
&\lesssim& \delta^{-(\rho+1+\epsilon)} (\sum_{j=0}^n (1+j)^{-2\rho-1-2\epsilon} |a_j|^2)^{1/2}.
\end{eqnarray*}
Since $\E (\sum_{j=0}^n (1+j)^{-2\rho -1-2\epsilon} |a_j|^2) = O(\sum_{j\ge 0} (1+j)^{-1-2\epsilon}) = O(1)$, we obtain
\begin{eqnarray*}
\E [\sup_{|z|\in I(\delta)+(-\delta/2,\delta/2)} |p_n(z)|^2] \lesssim \delta^{-2(\rho +1+\epsilon)}.
\end{eqnarray*}
The desired probabilistic estimate for $\log |p_n|$ then follows immediately. 

Now, the proof of the claimed probabilistic estimate for $\log |q_n|$ is similar. For convenience of notation, let $M_{\epsilon}=(\sum_{j=0}^n  (n+1-j)^{-2\rho}(1+j)^{-1-2\epsilon} |a_{n-j}|^2)^{1/2} $. Using Cauchy Schwarz and Lemma~\ref{l.elementary} we have, for $\epsilon>0$ and $w\in I(\delta)+ (-c\delta, c\delta)$:
\begin{eqnarray*}
|q_n(w)|  
&\lesssim& (n+1)^{-\rho} \Big(\sum_{j=0}^n (n+1-j)^{2\rho}(1+j)^{1+2\epsilon} |w|^{2j} \Big)^{1/2} M_\epsilon \\
&\lesssim&   \delta^{-(1+\epsilon)}  M_\epsilon.
\end{eqnarray*}
Again, $\E [M_\epsilon^2] \lesssim \sum_{j=0}^n  (1+j)^{-1-2\epsilon}   = O(1)$,  and the desired estimate follows immediately.
\endproof

\subsubsection{Proof of Theorem~\ref{t.locallogint}}

We will only show the proof for the claimed estimate for $\log |p_n|$, and the same argument works for $\log |q_n|$. Fix $z\in I(\delta)+(-c\delta, c\delta)$.  Let $C_1>0$ be big enough so that   Corollary~\ref{c.logpn} holds.

Thanks  Corollary~\ref{c.logpn}, we may assume that 
\begin{eqnarray}\label{e.lowerlogpn}
\log |p_n(z)| \ge -C_2|\log \delta|
\end{eqnarray}
for some $C_2>0$ large. Let $c'' \in (c', 1-c)$. Then for $w\in B(z,c''\delta)$ we have $|w|\in I(\delta)+(-(c+c'')\delta, (c+c'')\delta)$, so thanks to Lemma~\ref{l.unif-upper}, it holds with probability $1-O(\delta^{\alpha_0})$ that
\begin{eqnarray}\label{e.upperlogpn}
\sup_{w\in B(z,  c''\delta)} \log |p_n(w)| \le  C_3|\log \delta|
\end{eqnarray}
for  $C_3>0$ large.  

Below,  we will  condition on the event where \eqref{e.lowerlogpn} and \eqref{e.upperlogpn} hold,  on which we will show that
$$(\delta^{-2}\int_{B(z,c'\delta)} |\log |p_n(w)||^p dw)^{1/p} \lesssim  p |\log \delta|^2.$$

Now, the integrand $|\log |p_n||$ will blowup near the zeros of $p_n$, however only logarithmically. The above assumptions on $\log |p_n|$ will ensure that there are not many such zeros near $z$, and the main part of the argument is to control the zero-free part of $p_n$ using properties of subharmonic functions.

More specifically, let $\ell:=N_{p_n}(B(z,c''\delta))$ be the number of zeros of $p_n$ in $B(z, c'\delta)$.   As a consequence of Jensen's formula, we have
$$\ell   \lesssim_{c',c''} (\sup_{w\in B(z,c''\delta)} \log|p_n(w)| - \log |p_n(z)|) \lesssim |\log \delta|,$$
Now, let $u_1,\dots, u_{\ell}$ be the zeros of $p_n$ in $B(z, c'\delta)$. Let $Q_n(w) = p_n(w)/((w-u_1)\dots (w - u_\ell))$, this is a (random) polynomial having no zeros inside $B(z,c'\delta)$, we view $Q_n$ as the zero-free part of $p_n$.  It follows that, for any $p\ge 1$,
\begin{eqnarray*}
&&(\delta^{-2}\int_{B(z,c'\delta)} |\log |p_n(w)||^p)^{1/p}  \\
&\le& (\delta^{-2}\int_{B(z,c'\delta)} |\log |Q_n(w)||^p)^{1/p} + \sum_{i=1}^{\ell} (\delta^{-2}\int_{B(z,c'\delta)} |\log |w-u_i||^p)^{1/p} \\
&\lesssim&(\delta^{-2}\int_{B(z,c'\delta)} |\log |Q_n(w)||^p)^{1/p} + \ell p |\log \delta|.
\end{eqnarray*}
Since $\ell=O(|\log \delta|)$, it  remains to bound the integral involving   $Q_n$. In fact,  we will show that  $|\log |Q_n(w)|| = O(|\log \delta|^2)$ uniformly on $B(z,c'\delta)$, which is a stronger estimate. To see this, we first show that $\log |Q_n|$ satisfies   inequalities similar to \eqref{e.lowerlogpn} and \eqref{e.upperlogpn}. Indeed, note that  $\log |Q_n(w)|: B(0,c''\delta) \to \R \cup \{-\infty\}$ is a subharmonic function, and by the maximum principle it achieves its maximum on the boundary. It follows that
\begin{eqnarray*}
\sup_{w\in B(z,c''\delta)} \log |Q_n(w)| &\le& \sup_{w: \ |w-z|=c''\delta} \log |Q_n(w)|\\
&\le& \sup_{w:\ |w-z|=c''\delta}\log |p_n(w)| + \sup_{w: \ |w-z|=c''\delta} \sum_{i=1}^{\ell} \log \frac{1}{|w-u_i|} \\
&\lesssim& |\log \delta| + \ell |\log(\delta)|   \lesssim |\log \delta|^2.
\end{eqnarray*}
On the other hand, since $|z-u_i|\le c'\delta\le 1$ for all $i=1,\dots, \ell$, we also have
$$\log |Q_n(z)| = \log |p_n(z)| + \sum_{i=1}^\ell \log \frac{1}{|z-u_i|} \ge \log |p_n(z)| \ge  -C_2 |\log \delta|.$$
Thus we have verified that $Q_n$ satisfies inequalities similar to \eqref{e.lowerlogpn} and \eqref{e.upperlogpn}. Now, let $h(w):=C|\log \delta|^2 - \log|Q_n(w)|$ for a big constant $C$ such that  $h$ is nonnegative (and harmonic) on $B(z,c''\delta)$. Note that
$$0\le h(z) \le C|\log\delta|^2 + C_2|\log \delta| =O(|\log \delta|^2).$$
Using Harnack's inequality, for any $w\in B(z,c' \delta)$ we have
$$0\le h(w) \le \frac{c''\delta +c'\delta}{c''\delta -c'\delta}h(z)   = O(h(z))  = O(|\log \delta|^2).$$
It follows that  $|\log |Q_n(w)|| \le O( |\log \delta|^2) + |h(w)|  =O(|\log\delta|^2)$ for any $w\in B(z,c'\delta)$, as desired.  

\section{Counting local real roots}\label{s.localcount}

In this section, we will use the log integrability estimates  and the anti concentration estimates from previous sections   to establish several estimates for the local number of real roots for $p_n$.


For each $U\subset \C$ and any function $f$ analytic on a neighborhood of $U$, let $N_{f}(U)$  denote the number of roots of $f$ inside $U$.  

In this section, we assume that the coefficients of $p_n$  satisfy Condition~\ref{cond.polyvar}, and all implicit constants may depend on the implicit constants in Condition~\ref{cond.polyvar}. 

\begin{theorem}\label{t.localcount}  Let $0\le c,c'<1$ be such that $c+c'<1$. Then there are constants  $C_1,C_2,C_3>0$ such that the following holds: for any $\frac 1 n \lesssim \delta \le \frac 1 {C_1}$ and any $|z|\in I(\delta)+(- c\delta, c\delta)$ and any $M>0$ and  any event $E$ we have
\begin{eqnarray}\label{e.localcount}
\E [1_{E}N_{p_n}(B(z, c'\delta))^k] &\lesssim_{k,M}& \delta^M +  |\log \delta|^{C_2 k} \P(E) .
\end{eqnarray}
The analogous estimate also holds for $N_{q_n}=N_{p^*_n}$. Furthermore, for $\delta \ge  C_3 \log n/n$ we could take $C_2=1$.
\end{theorem}

It follows from Theorem~\ref{t.localcount}  that  the number of  roots of $p_n$ and $p^*_n$ on $I_{\R}(\delta)$ are at most logarithmic  away from $O(1)$. We state a useful corollary, when $E^c=\emptyset$.

\begin{corollary}\label{c.localcount} Let $0\le c,c'<1$ be such that $c+c'<1$. Then there are  constant  s $C_1,C_2,C_3>0$ such that for any $\frac 1 n \lesssim \delta \le \frac 1 {C_1}$ and any $|z|\in I(\delta)+(- c\delta, c\delta)$ we have
\begin{eqnarray}\label{e.localcount-cor}
\E [N_{p_n}(B(z, c'\delta))^k] \lesssim_{k}  |\log \delta|^{C_2k} .
\end{eqnarray}
Furthermore, for $\delta \ge C_3 \log n/n$  we could take $C_2=1$.
\end{corollary}

We  will divide the proof of Theorem~\ref{t.localcount}  into two cases, depending on whether $\delta$ is small or large. More specifically, we will consider first $\delta \ge  C_3 \log n/n$ for some sufficiently large  constant $C_3$, this is the large scale setting. Then we will consider the case when $\frac 1 n\lesssim \delta \lesssim \log n/n$ and refer to this as the small scale setting.  


\subsection{Larger scales}
We will use the following sublevel set estimate.

\begin{lemma}\label{l.sublevel-N}Let $0\le c,c'<1$ be such that $c+c'<1$. Let $C>0$ be sufficiently large. Let $\delta \in [\frac {C\log n} {n}, \frac 1 C]$  and assume that $|z|\in I(\delta)+(-c\delta,  c\delta)$.  Then  uniformly over $\lambda>C|\log  \delta|$ we have
\begin{eqnarray*}
\P(N_{p_n}(B(z, c'\delta))> \lambda) &\lesssim&  e^{-\lambda/C} + e^{-n\delta/C}, \\
\P(N_{q_n}(B(z, c'\delta))> \lambda) &\lesssim&  e^{-\lambda/C} + e^{-n\delta/C}.
\end{eqnarray*}
\end{lemma}

Let $C_3$ be large compared to the constant $C$ from Lemma~\ref{l.sublevel-N}. Using Lemma~\ref{l.sublevel-N} , we will prove \eqref{e.localcount}  for  $\delta>  C_3 \log n/n$. We will only show the details for $N_{p_n}$, the same argument could be applied to $N_{q_n}$. Now,  for brevity let $N=N_{p_n}(B(z, c'\delta))$ and $F= \{N \ge C_3|\log \delta|\}$. Since $N \le n$ trivially, we obtain
\begin{eqnarray*}
\E (N^k 1_{F}) &=& k\int_{t>0} t^k\P(N 1_F > t) \frac{dt}t\\
&\lesssim&  \int_{t\lesssim |\log \delta|} t^k \P(F) \frac{dt}t + \int_{|\log \delta|  \lesssim t \lesssim  n} t^k [e^{-n\delta/C} + e^{-t/C}]\frac{dt}t\\
&\lesssim& |\log \delta|^k \P(F) + n^k e^{-n\delta/C}  +   \int_{t\gtrsim   |\log \delta|} t^{k-1}e^{-t/C}dt \\
&\lesssim& |\log \delta|^k  (\delta^{C_3/C} + n^{-C_3/C}) + n^{k-C_3/C} + \delta^{C_3/(4C)} \quad  \lesssim_M \quad \delta^{M}
\end{eqnarray*}
if  $C_3$ is sufficiently larger than $CM$. It follows that
\begin{eqnarray*}
\E (N^k 1_{E}) &\lesssim& |\log \delta|^k \P(E) + \E (N^k 1_{F})\\
&\lesssim& |\log \delta|^k \P(E) +\delta^{M}, \quad \text{as desired}.
\end{eqnarray*}

\subsubsection{Proof of Lemma~\ref{l.sublevel-N}} Let $c''\in (c', 1-c)$. Using Jensen's formula, we have
$$N_{p_n}(B(z,c'\delta)) \le \sup_{w\in B(z,c''\delta)} \log |p_n(w)| - \log |p_n(z)|$$
\begin{align*} \text{therefore} &&  \P(N_{p_n}(B(z,c'\delta)) > \lambda)  \le \end{align*}
\begin{eqnarray}\label{e.app-Jensen}
&\le& \P(\sup_{w\in B(z, c''\delta)}  \log |p_n(w)| > \lambda/2) + \P(|p_n(z)|\le e^{-\lambda/2}).
\end{eqnarray}

For the first term on the right hand side of \eqref{e.app-Jensen}, we apply Lemma~\ref{l.unif-upper} with $s=\lambda/2$ and note that  $e^{2s}$ is a lot larger than any given power of $(1/\delta)$.

For the second term on the right hand side of \eqref{e.app-Jensen}, we use Theorem~\ref{t.smallball1}  with $t=e^{-\lambda/2}$ and use the assumption that $\lambda \ge C\log(1/\delta)$  (where $C$ is very large) to get the desired estimate.

The proof for $N_{q_n}$ is entirely similar.

\subsection{Smaller scales}
We now consider the smaller (and more critical) range  $\frac 1 n \lesssim \delta  \lesssim \frac{\log n}n$. Here  we will use Theorem~\ref{t.nns-weak} (from Section~\ref{s.loginteg}) about the log integrability of $p_n$ and $q_n$, which shows that there is
an event $F$ with probability $\P(F)=O(e^{-cn})$ such that for any $q\ge 1$ we have
\begin{eqnarray}\label{e.logpn}
\E [1_{F^c} \int_{B(0,1)} |\log|p_n(w)||^q  dw] &\lesssim&  (C q)^{Cq}   (\log n)^{Cq},
\end{eqnarray}
where $C$ is sufficiently large, and the analogous estimate also holds for $\log |q_n|$. We will use these estimates to show the desired estimates for $\log |p_n|$  in this range of $\delta$, and the argument for $\log|q_n|$ is entirely similar. 
\begin{align*}\text{\quad To start, note that} && \E [N_{p_n}(B(z,c'\delta))^k1_{E \cap F}] \ \lesssim \ n^k\P(F) = O(n^k e^{-cn})\end{align*}
which is $O_M(\delta^M)$ for any $M>0$. Thus, we may assume without loss of generality that $E\subset F^c$.   For convenience, denote $U=B(z,c'\delta)$ and $\Omega :=B(z, c''\delta)$ where  $c''\in (c', 1-c)$. Let $\phi$ be a smooth function such that $1_{B(0,c')} \le \phi \le 1_{B(0, c'')}$ and let $\phi_\delta(.)=\phi(./\delta)$ denote the $L^\infty$-preserving dilation of $\phi$.  We now use Green's formula 
$$\phi(0) = -\frac 1 {2\pi} \int_{\C} (\log |w|)\Delta \phi(w)dw$$
where $dw$ is the Lebesgue measure on $\C$. It follows that 
\begin{eqnarray*}
N_{p_n}(U)  \ \le \  \ \sum_{\alpha\in Z} \phi_\delta(z-\alpha)    &=& -\frac 1 {2\pi}   \int_{\mathbb C}  (\log|p_n(w)|)\Delta \phi_{\delta}(z_j-w)dw
\end{eqnarray*}
\begin{align*} 
\text{therefore} && N_{p_n}(U)   \ \lesssim \ \delta^{-2}\int_{\Omega} |\log|p_n(w)||dw.
\end{align*}
Consequently, using H\"older's inequality, the following holds for any $p \ge 1$
\begin{eqnarray*}
\E [N_{p_n}(U)^k1_{E}] 
&\lesssim& \delta^{-2k} \E[1_{E} (\int_{\Omega} |\log|p_n||)^k] \\
&\lesssim&  \delta^{-2k} \P(E)^{1-1/p} \Big(\E [1_{E} (\int_{\Omega} |\log|p_n(w)||   dw)^{kp}]\Big)^{1/p}  \\
&\lesssim&   \delta^{-2k} \P(E)^{1-1/p}  \Big(|\Omega|^{kp-1}\E [1_{E} \int_{\Omega} |\log|p_n(w)||^{kp}  dw] \Big)^{1/p}.  
\end{eqnarray*}
Recall that $E\subset F^c$ and note that $\Omega\subset B(0,1)$ and $|\Omega|=O(\delta^2)$. Therefore, using \eqref{e.logpn} for $q=kp$, we obtain
\begin{eqnarray*}
\E [N_{p_n}(U)^k1_{E}]  &\lesssim& \delta^{-2k}   \P(E)^{1-1/p} |\Omega|^{k-\frac 1p}  (Ckp)^{Ck}  \log^{Ck} n\\
&\lesssim_k&\delta^{-2/p} \P(E)^{1-1/p}   p^{Ck}  (\log n)^{Ck} .
\end{eqnarray*}
 Choosing $p=\log n \approx \log(1/\delta)$, then $\delta^{-1/p} = O(1)$, therefore
$$\E [N_{p_n}(U)^k1_{E}] \quad \lesssim\quad    \P(E)^{1-\frac 1 p} |\log \delta|^{(2C)k}.$$ 
Now, if $\P(E) \le \delta^{2M}$, then it is clear that the last right hand side is $O(\delta^M)$. If $\P(E)\ge \delta^{2M}$ then it is clear that $\P(E)^{1/p} \gtrsim 1$, consequently 
$$\E [N_{p_n}(U)1_{E}] \quad \lesssim\quad \P(E) |\log\delta|^{2Ck}.$$
 This completes the proof of Theorem~\ref{t.localcount}.

\section{Lindeberg swapping and Tao-Vu replacement estimates}
Our goal in this section is to establish the following result, which is a simple extension of a replacement estimate in Tao--Vu \cite{tv2015} to non-centered polynomials.

\begin{lemma}\label{l.logswap} For any $C,\epsilon,C_0>0$ there is $0<C_1<\infty$ so that the following holds.

Let  $\xi_0,\dots,\xi_n, G_0,\dots, G_n$ be  independent  with   $\E |\xi_j|^{2+\epsilon}<C$ and $\E |G_j|^{2+\epsilon}<C$ such that $\xi_j$ and $G_j$ have matching moments up to second order, for at least $n-C$ indices $j$. Let $\delta\in (0,1)$,  $\alpha_1>0$, $w_1,\dots, w_m \in I(\delta)$, and $F: \R^m \to \C$  be such that

(i) $m \lesssim \delta^{-\alpha_1}$, and   $|\partial^\beta F| \le  \delta^{-\alpha_1}$ for  $|\beta| \le 3$;

(ii) for all $1\le i\le m$ and $0\le j\le n$ it holds that $ |c_j w_i^j|\lesssim \delta^{C_1 \alpha_1} (\sum_j |c_j w_i^j|^2)^{1/2}$. 
\begin{align*}
\text{Then} &&& |\E F(\log |p_{n,\xi}(w_1)|,\dots, \log|p_{n, \xi}(w_m)|) \\
&&  & \ -   \ \E F(\log |p_{n,G}(w_1)|,\dots, \log |p_{n,G}(w_m)|)|  \ \ \lesssim \ \ \delta^{C_0\alpha_1},
\end{align*}
where the implicit constant may depend on $\alpha_1,C_0,C_1,\epsilon$.
\end{lemma}

Without loss of generality we may assume that $G_j$ are Gaussian for all $j$.
Following \cite{tv2015}, we will prove Lemma~\ref{l.logswap} using  the   Lindeberg swapping argument. The following basic estimate captures some ideas of this argument.

\begin{lemma} [Basic Lindeberg swapping] \label{l.lindeberg}
Let $\epsilon, C>0$. Assume that $\xi_1,\dots, \xi_n$ and $\widetilde \xi_1,\dots,\widetilde \xi_n$ are independent such that  $\max_j \E |\xi_j|^{2+\epsilon} \le C$ and $\max_j \E |\widetilde \xi_j|^{2+\epsilon} \le C$.

Assume that $\xi_j$ and $\widetilde \xi_j$ have matching moments up to second order for any $j\not\in J_0$.  Here $J_0$ is a subset of $\{1,\dots,n\}$.

Assume that $H: \C^{n} \to \C$, such that, as a function on $\R^{2n}$, $H\in C^3$.
Then  for some $\widetilde  C$ finite positive depending on $C$ and $\epsilon$ we have:
$$|\E H( \xi_1, \dots, \xi_n)  - \E H(\widetilde \xi_1,\dots,  \widetilde \xi_n)|   \le \widetilde C \Big( M_2^{1-\epsilon} M_3^{\epsilon} +  |J_0|^{2/3}\|H\|_{sup}^{2/3}M_3^{1/3}\Big).$$
Here viewing as a function on $\R^{2n}$ we let $M_i:=\sum_{j=1}^{n} \sum_{m=0}^{i}\|(\partial_{2j-1})^{i-m}(\partial_{2j})^m H\|_{sup}.$

\end{lemma}


\proof Let $H_1=H(\xi_1,\dots, \xi_n)$, and  let $H_{j+1}$ be obtained from $H_{j}$ by swapping $\xi_j$ with $\widetilde \xi_j$. We then estimate the left hand side by $\sum_j |\E (H_j-H_{j-1})|$.

Let $j\not\in J_0$. We view $H(\dots, w_{j}, \dots)$ as a function of $Re(w_j)$ and $Im(w_j)$, denoted by $f_j$.   For convenience, let  $M_{j,i}:=\sum_{m=0}^{i}\|(\partial_1)^{i-m}(\partial_2)^m f_j\|_{sup}$.

We consider approximation of $f_j(x,y)$ using Taylor expansion around $(0,0)$ up to second order terms.  By simple interpolation, the error term in this approximation is bounded above $O(\max(|x|^{2+\epsilon},|y|^{2+\epsilon}) M_{j,3}^{\epsilon} M_{j,2}^{1-\epsilon})$.
Since $\xi_j$ and $\widetilde \xi_j$ are independent from the others and have matching moments up to second order and since $\E|\xi_j|^{2+\epsilon}\le C$, $\E |\widetilde \xi_j|^{2+\epsilon} \le C$, it follows from direct examination that
$$\E[H_{j+1}-H_j ] =  O(M_{j,3}^{\epsilon} M_{j,2}^{1-\epsilon}).$$
Summing these estimates over $j\not\in J_0$ and using H\"older's inequality, we obtain
$$\sum_{j\not\in J_0} |\E[H_{j+1}-H_j ]|  \ \lesssim \   (\sum_j M_{j,3})^{\epsilon}  (\sum_j M_{j,2})^{1-\epsilon} \ = \ M_3^\epsilon M_2^{1-\epsilon} .$$

Now, let $j\in  J_0$. Again we view $H$ as a function $f_j$ of   $Re(w_j)$ and $Im(w_j)$ and approximate it by Taylor expansion around $(0,0)$ up to first order terms. We similarly obtain $ |\E[H_{j+1}-H_j ]| \lesssim   M_{j,1} (\E |\xi_j|+\E|\widetilde \xi_j|)= O(M_{j,1})$. Using Kolmogorov's inequality \cite{k1939} and a simple application of H\"older's inequality we obtain
\begin{eqnarray*}
\sum_{j\in J_0} |\E[H_{j+1}-H_j ]| &\lesssim&  \sum_{j\in J_0} M_{j,1} \lesssim \sum_{j\in J_0} M_{j,0}^{2/3} M_{j,3}^{1/3} \lesssim  |J_0|^{2/3}\|H\|_{sup}^{2/3} M_3^{1/3}.
\end{eqnarray*}
\endproof

We now prove Lemma~\ref{l.logswap}. Let $\sigma(z)=\sqrt{Var[p_{n,\xi}(z)]} =  (\sum_{0\le j \le n} |c_j z^j|^2)^{1/2}$. Let $\widetilde F:\R^m \to \C$ be defined by
$\widetilde F(u_1,\dots, u_m)=F(u_1 + \log \sigma(w_1),\dots, u_m + \log \sigma(w_m))$.
Then we also have $|\partial^\alpha \widetilde F| \lesssim \delta^{-\alpha_1}$ for all partial derivatives of order $|\alpha|\le 3$.  

Let $M=C_2\log (1/\delta)$ for some large constant $C_2>0$ to be chosen later. 

We perform a decomposition of $\widetilde F = F_1 + F_2$ where  $F_1= \phi \widetilde F$ and $F_2=  (1-\phi)\widetilde F$,
where $\phi$ is constructed below. Then $\phi:\R^m\to \R$ is a smooth function supported on $\{(x_1,\dots, x_m)\in \R^m: \min x_j \ge -(M+1)\}$ and equals $1$ on $\{(x_1,\dots, x_m)\in \R^m: \min x_j \ge -M\}$, such that $\|\partial^\alpha \phi\|_\infty \lesssim m^{|\alpha|}$ for any multi-index $\alpha$. 

We plan to apply Lemma~\ref{l.lindeberg} to 
$$H(\xi_0,\dots, \xi_n)=F_1(\log f(w_1),\dots, \log f(w_m)), \quad f(z):=|p_n(z)|/\sigma(z),$$ 
Now, $|\partial^{\alpha}F_1|  \ \lesssim \ m^3 \delta^{-\alpha_1} \  \lesssim \ \delta^{-4\alpha_1}$ for  $|\alpha|\le 3$. Via explicit computations,
\begin{eqnarray*}
\frac{\partial}{\partial Re(\xi_k)} \log |f(z)| \ = \ Re(\frac{c_k z^k}{p_n(z)}), &&
\frac{\partial}{\partial Im(\xi_k)} \log |f(z)| \  =\  -Im(\frac{c_k z^k}{p_n(z)}).
\end{eqnarray*}
Now, on the support of $F_1$ we have $\displaystyle |\frac{c_k w_j^k}{p_n(w_j)}|  \lesssim e^M  \frac{|c_k w_j^k|}{\sigma(w_j)}$. 
Thus, for $x,y\in\{Re(\xi_k), Im(\xi_k)\}$ we have
\begin{eqnarray*}
 |(\frac{\partial}{\partial x})(\frac{\partial}{\partial y})   H|  &\lesssim& \sum_{\ell,j=1}^m |(\partial_{\ell} \partial_j)F_1 |    |\frac {c_{k} w_j^{k}}{p_n(w_j)}| |\frac {c_{k} w_\ell^{k}}{p_n(w_\ell)}| +  \sum_{j=1}^m |\partial_j F_1|   \frac{|c_kw_j^k|^2}{|p_n(w_j)|^2} \\
 &\lesssim& e^{2M} \delta^{-4\alpha_1} \sum_{\ell,j=1}^m   |\frac {c_{k} w_j^{k}}{\sigma(w_j)}| |\frac {c_{k} w_\ell^{k}}{\sigma(w_\ell)}|  .
\end{eqnarray*}
Summing over $k$ and  using Cauchy Schwartz, we obtain
\begin{eqnarray*}
M_2 
&\lesssim& e^{2M}\delta^{-4\alpha_1}  \sum_{\ell,j=1}^m (\sum_{k=1}^n \frac{|c_k w_j^k|^2}{\sigma(w_j)^2})^{1/2} (\sum_{k=1}^n \frac{|c_k w_\ell^k|^2}{\sigma(w_\ell)^2})^{1/2} \quad \lesssim \quad e^{2M} \delta^{-6\alpha_1}.
\end{eqnarray*}

Similarly, we estimate the third partial derivatives for $H$ and use these estimates to bound $M_3$. Here we will arrive at trilinear sums, so using the assumption $|c_j w_k^j/\sigma(w_k)| = O(\delta^{C_1\alpha_1})$ we eventually obtain
\begin{eqnarray*}
M_3 &\lesssim& \delta^{-4\alpha_1} e^{3M} \sum_{\ell,j,h} \sum_{k} |\frac {c_{k} w_j^{k}}{\sigma(w_j)}| |\frac {c_{k} w_\ell^{k}}{\sigma(w_\ell)}| |\frac {c_{k} w_h^{k}}{\sigma(w_h)}|\ \lesssim \   e^{3M}\delta^{(C_1-7)\alpha_1} .
\end{eqnarray*}
Now, we may assume   $\epsilon\le 1$.  Via Lemma~\ref{l.lindeberg}, we have the generous bound
\begin{eqnarray*}
 |\E F_1(\log f(w_1),\dots) -  \E F_1(\log f_G(w_1),\dots) |  
&\lesssim&    e^{3M} \delta^{(C_1\epsilon-11)\alpha_1}.
\end{eqnarray*}

We now reset $H(\xi_0,\dots, \xi_n) := (1-\phi)(\log f(w_1),\dots, \log f(w_m))$. The partial derivatives of $(1-\phi)$ are  $O(1)$ and are supported in  $\min(\log f(w_1) ,\dots, \log f(w_m)) \ge -M-1$. Consequently, via   the same consideration as before, we obtain 
\begin{eqnarray*}
&&|\E \Big[F_2(\log  f(w_1),\dots, \log f(w_m))\Big]|  \quad \lesssim \quad \E [\delta^{-\alpha_1}H(\xi_0,\dots, \xi_n)] \\
&\lesssim& |\E[\delta^{-\alpha_1}H(G_0,\dots, G_n)]| + O(e^{3M}\delta^{(C_1\epsilon-11)\alpha_1}) \\
&\lesssim& \delta^{-\alpha_1}  \sum_{j=1}^m \P(\frac{|p_{n,G}(w_j)|}{\sqrt{Var[p_{n,G}(w_j)]}} < e^{-M}) + O(e^{3M}\delta^{(C_1\epsilon-11)\alpha_1}) \\
&\lesssim& \delta^{-2\alpha_1}  e^{-M} + O(e^{3M}\delta^{(C_1\epsilon-11)\alpha_1}) ,
\end{eqnarray*}
here we have used the fact that $p_{n,G}(w_j)$ is  Gaussian and $m=O(\delta^{-\alpha_1})$. 
Collecting estimates,  we obtain
$$|\E F(\log |p_{n,\xi}(w_1)|,\dots) - \E F(\log |p_{n,G}(w_1)|,\dots)| \quad \lesssim \quad \delta^{-2\alpha_1} e^{-M} + e^{3M} \delta^{(C_1\epsilon-11)\alpha_1}.$$
We choose $M=C_2\alpha_1 \log(1/\delta)$ where $C_2 \ge C_0+2$, and $C_1>(11+3C_2+C_0)/\epsilon$, then it is clear that the last right hand side is $O(\delta^{C_0\alpha_1})$, as desired.
This completes the proof of Lemma~\ref{l.logswap}.

\section{Proof of universality  for complex correlation functions}

In this section we prove Theorem~\ref{t.complex-corr}.   Following the framework developed by Tao-Vu \cite{tv2015}, we will use the Monte Carlo sampling method (summarized in Lemma~\ref{l.montecarlo}) and the Lindeberg swapping argument (implemented in Lemma~\ref{l.logswap}).  Below, we will only prove the desired estimates for the correlation functions of $p_n$.   The same argument could be applied to $q_n =(n+1)^{-\rho} p^*_n$ to get the desired estimates for  $p^*_n$.

We will actually show the desired estimates when $\phi_\delta$ has the tensor structure, namely $\phi_\delta(w)=\phi_{1,\delta}(w_1)\dots \phi_{k,\delta}(w_k)$, furthermore for such $\phi_\delta$ we will only need to assume that  each $\phi_{j,\delta}$, viewed as a function on $\R^2$, is continuously differentiable up to second order and furthermore   $|\partial^\alpha \phi_{j,\delta}| \le O(\delta^{-|\alpha|})$ for $|\alpha|\le 2$.  The reduction from  general (i.e. non tensor) $\phi_\delta$ to this special set up could be carried out as follows: First, let $c'\in (c,1)$, and let $\phi_{j,\delta}$ be smooth and supported inside $B_{\C}(0,c'\delta)$ such that $\phi_{j,\delta}=1$ on $B_\C(0,c\delta)$, and as a function on $\R^2$ it is $C^2$ and  satisfies the derivative bound $|\partial^\alpha \phi_{j,\delta}| \le O(\delta^{-|\alpha|})$  up to order $2$.
We may write
\begin{eqnarray*}
\phi_\delta (w_1,\dots, w_k) &=& \phi_{1,\delta}(w_1)\dots \phi_k(w_k) \phi(w_1,\dots, w_k)\\
 &=& \phi_1(w_1)\dots \phi_k(w_k) \sum_{n=(n_1,\dots, n_k)\in \mathbb Z^k} c_{n}e^{i 4\pi \delta^{-1}n\cdot w}\\
 &=&\sum_{n=(n_1,\dots, n_k)\in \mathbb Z^k} c_{n_1,\dots, n_k}  (\phi_1(w_1) e^{4\pi i n_1 w_1/\delta})\dots (\phi_k(w_k) e^{4\pi i nkw_k/\delta})
\end{eqnarray*}
using the multiple Fourier series expansion of $\phi$ on the polydisk $B_{\C}(0,\delta)^k$. By standard stationary phase estimates, if $\phi_\delta$ is $C^{m}$ then $|c_n| \lesssim_m (1+|n_1|+\dots +|n_k|)^{-m}$, while  $\partial^\alpha [\phi_j(w_j)e^{4\pi in_j w_j/\delta}] = O(\delta^{-|\alpha|}(1+|n_j|)^{|\alpha|})$, therefore if $m$ is large enough depending on $k$, say $m\ge 3k+2$, then we could write $\phi$ as a linear average of  tensor-type functions with the properties mentioned earlier.

Thus, we may now  assume that $\phi$ has the tensor structure. Let $z=(z_1,\dots,z_k)\in I(\delta)^k$ be fixed (no implicit constants  will depend on $z_j$'s). Recall that $Z$ denotes the multi-set of zeros of $p_n$. By definition,
$$\int_{\C^k} \phi_\delta(z-w) d\sigma(w)= \E \sum_{\alpha_1,\dots, \alpha_k\in Z} \phi_{1,\delta}(z_1-\alpha_1)\dots \phi_{k,\delta}(z_k-\alpha_k)$$
where the sum is over non repeated tuples of $k$ elements of the zero sets of $p_n$.
An application of the inclusion-exclusion formula will allow us to rewrite
the last right hand side
 as a linear combination of  terms, and each term is a product of finitely many sum of the following type 
 $$ X= \sum_{\alpha\in Z} \phi_{j,\delta,X}(z_j-\alpha),$$
 where $1\le j\le k$ is fixed and $\phi_{j,\delta,X}$ is a function supported in $B_{\C}(0, c\delta)$ such that, as a function on $\R^2$, it is  $C^2$ and its partial derivatives up to order $2$ are bounded accordingly.   
 
Consequently, it suffices to show that, for a sequence $X_{i_1}, \dots, X_{i_\ell}$ of the above type,
 $$|\E  X_{i_1}\dots X_{i_\ell} - \E   X_{G,i_1}\dots X_{G,i_\ell}| = O(\delta^c)$$
(uniform over all choices of $1\le \ell\le k$ and $1\le i_1 <\dots < i_\ell\le k$), for some $c>0$. Without loss of generality, we may  assume that $\ell=k$ and $i_1=1$,... ,$i_k=k$, and for brevity we will omit the dependence on $X_j$ in the notation and simply write $X_j=\sum_{\alpha}\phi_{j,\delta}(z_j-\alpha)$ below. 

Let $\alpha_0>0$ be a sufficiently small constant that  may depend  on the underlying implicit constants   in Condition~\ref{cond.polyvar}. 
By a standard construction, we could find $\varphi:\C^k \to \C$ such that $\phi$ supported on $B(0,2 \delta^{-\alpha_0})$ and $\varphi(w_1,\dots, w_k)=w_1\dots w_k$ on $B(0,\delta^{-\alpha_0})$, furthermore  $|\varphi(w_1,\dots,w_k)|\le |w_1\dots w_k|$ for any $w_1,\dots, w_k$, and   (as a function on $\R^{2k}$)  $\varphi$ will be in $C^2$ with $|\partial^\alpha \varphi(w)| \lesssim \delta^{-k\alpha_0}$ for any (partial) derivatives of order up to $2$. 

Let $C>0$ is sufficiently large and let $\frac 1 n \lesssim \delta \le \frac 1 C$.  We first use Theorem~\ref{t.locallogint} and Lemma~\ref{l.unif-upper} to conclude that for any $0<c'<1/2$ there is an event $E=E(\delta,\alpha_0,z_1,\dots, z_k)$ with probability $\P(E) =  O_{c',\alpha_0} (\delta^{c'})$ such that on $T=E^c$ the following holds for each $j=1,2,\dots, k$:
\begin{eqnarray*}
\sup_{w: |w-z_j|\le c\delta}\log |p_n(w)|  &\lesssim& |\log \delta|, \\
\frac 1 {\delta^2}\int_{B(z_j, c\delta)} |\log |p_n(w)||^2 dw &\lesssim& |\log \delta|^4.
\end{eqnarray*}

We now use Green's formula, which says that the following holds for any $\phi$ compactly supported in $C^2(\R^2)$
$$\phi(0) = -\frac 1 {2\pi} \int_{\C} (\log |w|)\Delta \phi(w)dw$$
where $dw$ is the Lebesgue measure. It follows that, for each $1\le j\le k$, we have
\begin{eqnarray}
\nonumber X_j &=&  \sum_{\alpha\in Z}(-\frac 1 {2\pi}) \int_{\mathbb C} \log |w-\alpha| \Delta \phi_{j,\delta}(z_j-w)dw\\
\label{e.Xjinteg} &=& -\frac 1 {2\pi}   \int_{\mathbb C}  (\log|p_n(w)|)\Delta \phi_{j,\delta}(z_j-w)dw.
\end{eqnarray}
Thus, using H\"older's inequality and using the above properties of $T$, we obtain $|X_j| \lesssim |\log \delta|^2$ on the event $T$. By ensuring that $\delta<1/C$ for $C$ sufficiently large, it follows that $|X_j|<\delta^{-\alpha_0}$ on the event $T$. Now, outside $T$ we still have $|\phi(X_1,
\dots,X_k)|\le |X_1\dots X_k|$, therefore
\begin{eqnarray}\label{e.preswap}
\E X_1\dots X_k &=&   \E \varphi(X_1,\dots, X_k) + O(\max_j \E[|X_{j}|^k1_{E}]).
\end{eqnarray}

We now use Monte Carlo sampling to approximate the integral form \eqref{e.Xjinteg} of $X_j$ with a discrete sum. 

\begin{lemma}[Monte Carlo sampling]\label{l.montecarlo} Let $(X,\mu)$ be a probability space and let $f \in L^2(X,\mu)$. Assume that $w_1,\dots, w_m$ are drawn independently from $X$ using the distribution $\mu$. Then for $S=\frac 1 m(f(w_1)+\dots + f(w_m))$ we have $\E S =  \int_X f d\mu$ and
$$\P (|S - \E S | \ge \lambda) \le \frac 4 {m \lambda^2}\int_X |f|^2 d\mu.$$
\end{lemma}

Now, $\Delta \phi_{j,\delta}$ is supported inside $B(0, c\delta)$ and is bounded above by $O(\delta^{-2})$. 

Let $w_{j,i}$ be uniformly chosen from $B(0,c\delta)$ (independent of each other and of the coefficients of $p_n$), here $1\le i\le m$ and $1\le j\le k$. 
Using \eqref{e.Xjinteg} and  Lemma~\ref{l.montecarlo}, it follows that
\begin{eqnarray*}
\P(|X_j-\frac 1 m \sum_{i=1}^m a_{j,i}\log|p_n(w_{j,i})||>\lambda) &\lesssim& m^{-1} \lambda^{-2}  \delta^{-2} \int_{B(z_j, c\delta)}  |\log  |p_n(w)||^2 dw,
\end{eqnarray*}
where $a_{j,i}= -\frac 1 {2} c^2\delta^2 \Delta \phi_{j,\delta}(z_j-w_{j,i})$. Note that $|a_{j,i}|=O(1)$.

Now, on the event $T$, the right hand side in the last display is   $O(m^{-1} \lambda^{-2}|\log \delta|^4)$.  Using the above estimate, we now show that all $X_j$'s could be replaced by   the corresponding averages at a total small cost:
\begin{claim}\label{cl.swapinsidephi}
Let $w=(w_{11},\dots, w_{1m},\dots, w_{k1},\dots, w_{km})$. Then
$$|\E  \varphi(X_1, .., X_k) - \E  \varphi(\frac 1 m \sum_{i=1}^m a_{1,i}\log|p_n(w_{1,i})|,.., \frac 1 m \sum_{i=1}^m a_{k,i}\log|p_n(w_{k,i})|)| \ =\  O(\delta^{\alpha_0}),$$
where the expectation is taken over   $w$ and  $\xi=(\xi_0,\dots, \xi_n)$.
\end{claim}

To see this, let $\lambda = \delta^{(k+1)\alpha_0}$.  Then on the product probability space generated by $\xi=(\xi_0,\dots, \xi_n)$ and $w_j=(w_{j,1},\dots, w_{j,m})$ it holds with probability $1-\P(T^c)- O_k(m^{-1}\delta^{-(2k+3)\alpha_0})$ that
$$\Big|X_j-\frac 1 m \sum_{i=1}^m a_{j,i}\log|p_n(w_{j,i})| \Big| \lesssim \delta^{(k+1)\alpha_0},$$
for all $j=1,\dots, k$. Now, letting   $m\approx \delta^{-(3k+4)\alpha_0}$ and choosing $\alpha_0$ sufficiently small (so that in particular $c>(k+1)\alpha_0)$), it follows that the following inequality holds with probability $1- O(\delta^{(k+1)\alpha_0})$:
$$| \varphi(X_1,\dots, X_k) -  \varphi(\frac 1 m \sum_{i=1}^m a_{1,i}\log|p_n(w_{1,i})|, \dots )| = O(\delta^{\alpha_0}).$$
(Here we've used the assumption that the first order partial derivatives of $\varphi$ is bounded above by $O(\delta^{k\alpha_0}))$.)
On the event that this estimate does not hold (which has probability $O(\delta^{(k+1)\alpha_0})$), we have the crude bound $O(\delta^{-k\alpha_0})$ for the left hand side of the above display, here we have used the assumption that $|\phi(w_1,\dots,w_k)|\le |w_1\dots w_k|$ and $\phi$ is supported on $B_{\C}(0,2\delta^{-\alpha_0})^k$. 
Collecting estimates,  the desired estimate of Claim~\ref{cl.swapinsidephi} follows immediately.

On the event $E$, we note that $X_j \lesssim N_{p_n}(B(z_j, c\delta))$ and similarly $X_{j,G} \lesssim N_{p_{n,G}}(B(z_j, c\delta))$. Consequently, using \eqref{e.preswap} and Claim~\ref{cl.swapinsidephi}  we obtain
\begin{eqnarray}
\label{e.preLindeberg}  |\E X_1\dots X_k - \E X_{1,G}\dots X_{k,G}|   = 
\end{eqnarray}
$$= |\E  \varphi(\frac 1 m \sum_{i=1}^m a_{1,i}\log|p_n(w_{1,i})|, \dots) -\E\varphi(\frac 1 m \sum_{i=1}^m a_{1,i}\log|p_{n,G}(w_{1,i})|, \dots) | + $$
$$  + O(\sum_j \E[ 1_E N_{p_n}(B(z_j, c\delta))^k]) + O(\sum_j  \E[ 1_E N_{p_{n,G}}(B(z_j, c\delta))^k]) + O(\delta^{\alpha_0}).$$

Using Theorem~\ref{t.localcount},   the two terms involving  $N_{p_n}(B(z_j, c\delta))$ and $N_{p_{n,G}}(B(z_j, c\delta))$ are bounded by $O(|\log \delta|^{Ck} \delta^{\alpha_0})$, which in turn is bounded by $O(\delta^{\alpha_0/2})$. 

Thus, it remains to bound  the first term on the right hand side of \eqref{e.preLindeberg}. Here we use Lindeberg swapping, or more precisely Lemma~\ref{l.logswap}. Below we  only discuss swapping of $\frac 1 m\sum_{i=1}^m \log |p_n(w_{1,i})|$ with its Gaussian analogue $\frac 1 m \sum_{i=1}^m a_{1,i}\log|p_{n,G}(w_{1,i})|$; the swapping of the other $k-1$ averages can be done similarly.   Now, by conditioning on other variables and treating them as parameters, we may let
$$F(u_1,\dots, u_m) = \phi(\dots, \frac 1 m (a_{1,i}u_1+\dots+ a_{1,m}u_m), \dots).$$
It remains to show that
$$\E F(\log|p_n(w_{1,1})|, \dots, \log |p_n(w_{1,m})|) - \E F(\log|p_{n,G}(w_{1,1})|, \dots, \log |p_{n,G}(w_{1,m})|)  $$
$$\lesssim \delta^{\alpha_0}.$$
We can check that $|\partial^\beta F| \lesssim \frac 1 {m^{|\beta|}} \delta^{-k\alpha_0}$ for any partial derivatives up to order $3$. Note that  $m\approx \delta^{-(3k+4)\alpha_0}$ by choice and $\alpha_0$ could be chosen arbitrarily small. Therefore, in order to show the estimate in the last display via Lemma~\ref{l.logswap}, it remains to show that  for some uniform constant $c>0$ (independent of $\alpha_0$) the following holds
$$|c_jw_i^j| \lesssim \delta^c \sqrt{Var[p_n(w_i)]}$$ 
for any $1\le i\le m$ and any $0\le j\le n$. To see this, note that $1-|w_i|\approx \delta$ and $c_j$'s satisfy Condition~\ref{cond.polyvar}, therefore
$$\sqrt{Var[p_n(w_i)]} \gtrsim (\sum_{j} j^{2\rho} |w_i|^{2j})^{1/2} \gtrsim \sqrt{(1-|w_i|^2)^{-2\rho-1}} \gtrsim \delta^{-\rho-1/2},$$
while $|c_j w_i^j| \lesssim (1+ j)^\rho (1-\delta)^j$. Via examination of the function $x^\rho (1-\delta)^x$ over $x\in [0,\infty)$, we could  show that $|c_j w_i^j|/\sqrt{Var[p_n(w_i)]} \lesssim \delta^{\rho+\frac 1 2} + \delta^{1/2}$, thus we could take any $0<c\le \min(\rho+\frac 12, \frac 1 2)$. (Recall the assumption that $\rho>-1/2$).


%
%

\section{Counting local non-real roots}

In this section, we will prove several  estimates for the local number of  non-real roots of $p_n$ near the real line. These estimates play an essential role in the next section, where the proof of   Theorem~\ref{t.corr} will be presented.
Recall that we write $p_n=m_n+r_n$  where $m_n(z)=\sum_j b_j z^j$ is the deterministic component and $r_n=\sum_j c_j \xi_j z^j$ is the random component. We divide the analysis into two scenarios.

{\bf Scenario 1:     $m_n$ is ``small'' compared to   $r_n$.}  
This scenario generalizes  the special case $m_n=0$ considered in in \cite{dnv2017}, where it was shown that with high probability $r_n$ has no non-real local root. Here we will show that a similar conclusion holds even with the addition of a ``small" deterministic component $m_n$.

\begin{lemma}\label{l.nearR} Let $\epsilon_0>0$ be sufficiently small   and let $c\in [0,1)$. Then for $C=C(\epsilon_0,c)>0$ sufficiently large the following holds for any $\frac1 n \lesssim \delta \le \frac 1 C$ and $\eta:=\delta^{1+\epsilon_0}$ and any $x\in  I_{\R}(\delta)+ (-c\delta, c\delta)$. 

(i)  Assume that on $B(x,2\eta)$ we have $|m''_n| \ \lesssim \    \sqrt{Var[r''_n]}$.
\begin{align*}
\text{Then for any $\kappa<2$ we have} && \P(N_{p_n}(B(x, \eta))\ge 2) = O_{\epsilon_0,\kappa}  ((\eta/\delta)^\kappa) .
\end{align*}
(ii)   Assume that on $B(x,2\eta)$ we have $|{m^*}''_n|  \ \lesssim \   \sqrt{Var[{r_n^*}'']}$.
\begin{align*}
\text{Then for any $\kappa<2$ we have} && \P(N_{p^*_n}(B(x, \eta))\ge 2) = O_{\epsilon_0,\kappa}  ((\eta/\delta)^\kappa) .
\end{align*}
\end{lemma}
 
{\bf Scenario 2:    $m_n$  is ``large'' compared to $r_n$}.
Here we will show that  with high probability $p_n$ has no local roots  in a neighborhood of  the real line. 

 \begin{lemma}\label{l.nearR-mlarge} 
Let $\epsilon_0>0$ be  sufficiently small   and let $c\in [0,1)$. Let $\kappa>0$. Then for $C, C'>0$ sufficiently large   the following holds for any $\frac1 n \lesssim \delta \le \frac 1 C$ and $\eta :=\delta^{1+\epsilon_0}$ and any $x\in  I_{\R}(\delta)+(-c\delta,c\delta)$.

 (i) Assume that on $B(x,2\eta)$ we have $|m_n|  \  > \ C' |\log  \delta|^{1/2}  \sqrt{Var[r_n]}$.
 \begin{align*}
 \text{Then}  && \P(N_{p_n}(B(x, \eta))\ge 1) = O  ((\eta/\delta)^\kappa).
 \end{align*}

(ii) Assume that on $B(x,2\eta)$ we have $|{m^*}_n|  \ > \ C' |\log \delta|^{1/2}  \sqrt{Var[r^*_n]}$.
\begin{align*}
\text{Then} &&
\P(N_{p^*_n}(B(x,\eta)) \ge 1) =  O_{\epsilon_0,C_0,\kappa}  ((\eta/\delta)^\kappa) .
\end{align*}
\end{lemma}

\subsection{Proof of Lemma~\ref{l.nearR}}
\subsubsection{Proof of Lemma~\ref{l.nearR}, part (i)}  Here we prove part (i) and we will discuss the modifications for part (ii) later. For convenience, let
$$X=N_{p_{n}}(B(x,  \eta)), \ \ X_G=N_{p_{n,G}}(B(x,  2\eta)).$$

\underline{Step 1. Reduction to Gaussian:}
We'll  use  Theorem~\ref{t.complex-corr} in this step. Let $\widetilde c\in (c,1)$.

Let $\eta_1,\dots, $ be an enumeration of the (complex) roots  of $p_n$ and let $\eta_{1,G}, \dots$ be an enumeration of the (complex) roots of $p_{n,G}$, both enumerated with multiplicity. 

Let $\epsilon_1>0$ be   small  to be chosen later.  Let $\varphi:\C \to [0,1]$ be smooth supported on $B(0,2)$ such that $\varphi(z)=1$ if $|z|\le 1$.
We have
$$ \P(X \ge 2) \le \E \sum_{i\ne j} \varphi(\frac{\eta_i-x}{\eta}) \varphi(\frac{\eta_j -x}{\eta}).$$
We now discuss the set up required to apply Theorem~\ref{t.complex-corr}. Since $x\in I_{\R}(\delta)+(-c\delta,c\delta)$, we may write $x=x_0+ \alpha$ where $x_0\in I_{\R}(\delta)$ and $|\alpha|\le c\delta$. We then let
$$\phi_{\delta}(z,w):=\delta^{L\epsilon_0}\varphi(\frac {z-\alpha}{\eta})\varphi(\frac {w-\alpha}{\eta})$$
which is defined on $\C^2$, and here $L=O(1)$ is a sufficiently large absolute constant (in particular independent of $\epsilon_0$) so that all required derivative bounds (from Theorem~\ref{t.complex-corr}) for $\phi_\delta$ are satisfied. Now, $supp(\phi_\delta)\subset B_{\C}(\alpha,2\eta)^2 \subset B_{\C}(0, \widetilde c \delta)^2$ if we require $\delta<1/C$ with $C>0$ sufficiently large depending on $\epsilon_0$,$c$, and $\widetilde c$.   It then follows from Theorem~\ref{t.complex-corr} (and the definition of correlation functions) that for some $\alpha_0>0$ (independent of $L,\epsilon_0$) the following holds:
\begin{eqnarray*}
\E \sum_{i\ne j}  \phi_\delta(\eta_i - x_0,  \eta_j-x_0)  &=&  \E \sum_{i\ne j}  \phi_\delta(\eta_{i,G} - x_0,  \eta_{j,G} -x_0) + O(\delta^\alpha_0).
\end{eqnarray*}
Unraveling the notation, we obtain
\begin{eqnarray*}
\E \sum_{i\ne j} \varphi(\frac{\eta_i-x}{\eta}) \varphi(\frac{\eta_j -x}{\eta}) &\le&  \E  \sum_{i\ne j} \varphi(\frac{\eta_{i,G}- x}{\eta}) \varphi(\frac{\eta_{j,G} - x}{\eta}) + O(\delta^{\alpha_0} \delta^{-L\epsilon_0}) \\
 &\le&   \E  [X_G(X_G-1)] + O(\delta^{\alpha_0-L\epsilon_0})  \\
&\le&  \delta^{-2\epsilon_1}\P(X_G\ge 2) + \E  [X_G^2  1_{X_G > \delta^{-\epsilon_1}}] + O(\delta^{\alpha_0-L\epsilon_0}) .
\end{eqnarray*}
Using Corollary~\ref{c.localcount}  and observing that $X_G\le N_{p_{n,G}}(B(z,\delta/9))$, we have 
$$\P(X_G > \delta^{-\epsilon_1}) \lesssim_m  \delta^{m\epsilon_1} |\log\delta|^{O(m)}$$ 
for any $m\ge 1$, so by choosing $m$ large we have a bound of $O_{\epsilon_1,M}(\delta^M)$ for any $M>0$. Using Theorem~\ref{t.localcount}, it follows that
$$\E [X_G^2 1_{X_G >\delta^{-\epsilon_1}}]  \quad \lesssim_{M} \quad \delta^M +  |\log \delta|^{O(1)} \P(X_G > \delta^{-\epsilon_1}) \quad \lesssim \quad \delta^{M/2}.$$
Collecting estimates, we obtain
\begin{eqnarray*}
\P(X \ge 2) &\le& \delta^{-2\epsilon_1} \P(X_G\ge 2)  + O(\delta^{\alpha_0-2L\epsilon_0}) \\
&\le& \delta^{-2\epsilon_1}\P(X_G\ge 2) + O(\delta^{\kappa \epsilon_0})
\end{eqnarray*}
by choosing $\epsilon_0$  small.  So it remains to show that $P(X_G\ge 2) \lesssim \delta^{\kappa \epsilon_0+2\epsilon_1}$. Since $\epsilon_1$ could be chosen very small, it suffices to show that $P(X_G\ge 2) \lesssim \delta^{\kappa' \epsilon_0}$ for some $\kappa' \in (\kappa,2)$, which is essentially the Gaussian analogue of the desired estimate. 


\underline{Step 2. Proof for  Gaussian.}
We will show that, with high probability $p_{G,n}$ is close to its linear approximation at $x$, namely $\mathcal L(z) := p_{n,G}(x) + p'_{n,G}(x)(z-x)$. 
\begin{eqnarray}\label{e.linearapprox}
\P(\min_{z\in \partial B(x,2\eta)} |\mathcal L(z)| \le \max_{z\in \partial B(x,2\eta)} |\mathcal E(z)|)  = O(\delta^{\kappa \epsilon_0}).
\end{eqnarray}
Using Rouch\'e's theorem and linearity of $\mathcal L$, \eqref{e.linearapprox} implies the desired estimate. Now, to show \eqref{e.linearapprox}, we will prove two estimates.

\begin{claim}\label{cl.concentrationL} 

The following holds uniformly over $t>0$:
$$\P(\min_{z\in \partial B(x,2\eta)} |\mathcal L(z)| \le t\delta^2 \sup_{\xi\in B(x,2\eta)}\sqrt{Var[r''_n(\xi)]})  = O(t).$$
\end{claim}

\begin{claim}\label{cl.deviationE} For some $\alpha_0>0$ the following holds uniformly over   $t>0$:
$$\P(\max_{z\in \partial B(x, 2\eta)} |\mathcal E(z)| > t \eta^2  \sup_{\xi\in B(x,2\eta)}\sqrt{Var[r''_n(\xi)]})  =O(e^{-\alpha_0 t^2})$$
\end{claim}
The desired estimate \eqref{e.linearapprox}  then follows from choosing $t=(\eta/\delta)^{\kappa}$ in Claim~\ref{cl.concentrationL} and choosing $t=M|\log\delta|^{1/2}$ (with $M$ large) in  Claim~\ref{cl.deviationE}. Here we need $\kappa<2$.

\subsubsection{Proof of Claim~\ref{cl.concentrationL}}

Since $\mathcal L$ is linear with real coefficients and since $x\in \R$,   $\min_{z\in \partial B(x,2\eta)} |\mathcal L(z)|$ is achieved at $z=x-2\eta$ or $z=x+2\eta$. Consequently, for any $t>0$ we have
\begin{eqnarray*}
\P(\min_{z\in \partial B(x,2\eta)} |\mathcal L(z)| \le t) 
&\le&  \P(|\mathcal L(x+2\eta)|\le t) + \P(|\mathcal L(x-2\eta)|\le t) \\
&\lesssim& \frac{t}{\sqrt{Var[\mathcal L(x+2\eta)]}} + \frac{t}{\sqrt{Var[\mathcal L(x-2\eta)]}},
\end{eqnarray*}
here we have used the fact that $\mathcal L(x+2\eta)$ and $\mathcal L(x-2\eta)$ are Gaussian. Using Lemma~\ref{l.elementary} and Condition~\ref{cond.polyvar}, we have $Var[r_n(x)] \approx \delta^4 \sup_{\xi\in B(x,2\eta)}Var[r''_n(\xi)]$.
Therefore it remains to show that for any $s\in \{-2\eta, 2\eta\}$ we have
\begin{eqnarray}\label{e.varL>varp}
\sqrt{Var[\mathcal L(x + s)]}  &\gtrsim& \sqrt{Var[r_n(x)]}.
\end{eqnarray}
Now, since $\xi_j$ are independent, we have $\sqrt{Var[\mathcal L(x + s)]}  =  \|c_j (x^j + s j x^{j-1})_{j=0}^n\|_{l^2}$.

If $\delta \ge \frac 1 {10n}$ then by definition we have $  1- |x| \approx\delta$. Therefore,  using the triangle inequality and Lemma~\ref{l.elementary} and Condition~\ref{cond.polyvar} we obtain
\begin{eqnarray*}
\sqrt{Var[\mathcal L(x  + s)]} 
&\ge&  \|(c_j x^j)_{j=0}^n\|_{l^2} -|s| \|(c_j j x^{j-1})_{j=0}^n\|_{l^2} \\
&\ge& \sqrt{Var[r_n(x)]} - 2\eta \sqrt{Var[r'_n(x)]} .
\end{eqnarray*}
Using Lemma~\ref{l.elementary} and Condition~\ref{cond.polyvar}, it follows that $ Var[r'_n(x)] \approx \delta^{-2} Var[r_n(x)]$. Since $\eta \ll \delta$, the desired estimate \eqref{e.varL>varp} follows immediately.

Now, if $\frac 1 n \lesssim \delta <\frac1{10n}$ we have $|s|\le 2\eta<1/(2n)$. Therefore, uniformly over $0\le j \le n$ we have $|x^j + s j x^{j-1}|  \gtrsim  |x|^j$,  which  implies the desired estimate \eqref{e.varL>varp}.

\endproof

\subsubsection{Proof of Claim~\ref{cl.deviationE}.}

To estimate $\max_{z\in \partial B(x, 2\eta)} |\mathcal E(z)|$, we first estimate the mean and the variance of $\mathcal E(z)$. We will show that
\begin{eqnarray}
\label{e.meanError}|\E \mathcal E (z)|  &\lesssim& \eta^2 \sup_{\xi\in B(x,2\eta)}\sqrt{Var[r''_n(\xi)]} , \\ 
\label{e.varError} Var [\mathcal E(w)] &\lesssim& \eta^4 \sup_{\xi\in B(x,2\eta)} Var[r''_n(\xi)]    
\end{eqnarray}
uniformly over $z\in B(x, 2\eta)$ and $w\in  B(x,3\eta)$. 


For \eqref{e.varError}, let $w\in B(x,3\eta)$. By the mean value theorem, we have
\begin{eqnarray*}
Var[\mathcal E(w)] &=& \sum_{j=0}^n   |c_j|^2 |w^j - x^j - j(w-x)x^{j-1}|^2\\
&\lesssim& \eta^4 \sup_{\xi\in B(x,3\eta)} Var[r''_n(\xi)] .
\end{eqnarray*}
By ensuring $C=C(\epsilon_0,c)$ is large, for any $\xi \in B(x,3\eta)$ we have $\xi\in I(\delta)+(-c'\delta,c'\delta)$ for $c'=(1+c)/2<1$.  Using Lemma~\ref{l.elementary}, it follows that 
$$\sup_{\xi\in B(x,3\eta)} Var[r''_n(\xi)]  \lesssim \sup_{\xi\in B(x,2\eta)} Var[r''_n(\xi)], \quad \text{which implies \eqref{e.varError}}.$$

For \eqref{e.meanError},  again by the mean value theorem we have
\begin{eqnarray*}
|\E \mathcal E(z)| &=& |m_{n}(z) - m_{n}(x) - m'_{n}(x)(z-x)| \\
&\lesssim& \eta^2 \sup_{\xi\in B(x,2\eta)} |m''_{n}(\xi)|\\
 &\lesssim& \eta^2 \sup_{\xi\in B(x,2\eta)}\sqrt{Var[r''_n(\xi)]}  \quad \text{(by the given assumption)}
\end{eqnarray*}

 Now, we combine \eqref{e.meanError} and \eqref{e.varError} to prove Claim~\ref{cl.deviationE}. For convenience of notation, let $q(z):=\mathcal E(z) - \E \mathcal E(z)$. Without loss of generality, we may assume that $t$ is much larger than  the implicit constants in the last estimate for $|\E \mathcal E(z)|$ and in \eqref{e.meanError}. It follows from \eqref{e.meanError} and \eqref{e.varError} that  
\begin{eqnarray*}
&&\P(\max_{z\in \partial B(x, 2\eta)}|\mathcal E(z)| \ge  t\eta^2 \sup_{\xi\in B(x,2\eta)}\sqrt{Var[r''_n(\xi)]} ) \\
&\le& 
\P (\max_{z\in \partial B(x, 2\eta)} |q(z)| \ge (t/2)\eta^2 \sup_{\xi\in B(x,2\eta)}\sqrt{Var[r''_n(\xi)]} )\\
&\le& \P (\max_{z\in \partial B(x, 2\eta)} |q(z)| \gtrsim t   \sup_{\xi\in B(x,2\eta)}\sqrt{Var[\mathcal E(\xi)]} ).
\end{eqnarray*}
Using Cauchy's theorem, for $z\in \partial B(x, 2\eta)$ we have
\begin{eqnarray*}
|q(z)| 
&\lesssim&   \int_{\partial B(x,3\eta)} |q(w)|\frac{d|w|}{\eta} \\
&\lesssim&    \sup_{w\in B(x,3\eta)} \sqrt{Var[q(w)]} \int _{\partial B(x,3\eta)} \frac{|q(w)|}{\sqrt{Var[q(w)]}} \frac{d|w|}{\eta}
\end{eqnarray*}
where $d|w|$ is the arclength measure along the integration contour $\partial B(x,3\eta)$.
Note that $Var[q(z)]=Var[\mathcal E(z)]$. It follows that, for some $c>0$, we have

\begin{eqnarray*}
&&\P(\max_{z\in \partial B(x, 2\eta)} |q(z)| \ge t \sup_{\xi\in B(x,3\eta)}\sqrt{Var[\mathcal E(\xi)]})\\
&\lesssim&  e^{-ct^2} \E \exp(\int _{\partial B(x,3\eta)} \frac{|q(w)|}{2\sqrt{Var [q(w)]}}\frac{d|w|}{\eta})^2\\
&\lesssim& e^{-ct^2} \int _{\partial B(x,3\eta)} \E \exp(\frac{|q(w)|^2}{ 4Var [q(w)] })\frac{d|w|}{\eta} \quad \text{(by convexity)} \\
&\lesssim& e^{-ct^2}  \quad \text{ (since $\frac{q(w)}{\sqrt{Var[q(w)]}}$ is normalized Gaussian)}.
\end{eqnarray*}
\endproof

\subsubsection{Proof of Lemma~\ref{l.nearR}, part (ii)}

Our proof of part (ii) of Lemma~\ref{l.nearR}  is entirely similar to that of the proof of part (i), where the key ingredients is the fact that uniformly over $\xi\in B(x,3\eta)$ we have $\sqrt{Var[{r^*}^{(m)}_n(\xi)]} \approx_m  (1+n)^\rho \delta^{-(2m+1)/2}$ for any $m\ge 0$, which in turn is a consequence of Condition~\ref{cond.polyvar} and Lemma~\ref{l.elementary}.  

\endproof

\subsubsection{Proof of Lemma~\ref{l.nearR-mlarge}, part (i)}
We will proceed in a similar fashion as in the proof of Lemma~\ref{l.nearR}. The reduction to the Gaussian setting can be done similarly by using universality estimates for the 1-point correlation function of the complex zeros of $p_n$ from  Theorem~\ref{t.complex-corr} and estimates proved in Theorem~\ref{t.localcount} and Corollary~\ref{c.localcount}. 

We now discuss the proof for the Gaussian setting. The given assumption clearly implies that $m_n$ has no zero in $B(x,2\eta)$. Thus, using Rouch\'e's theorem it suffices to show that
$$P(\sup_{\xi \in \partial B(x,2\eta)} |r_{n,G}(\xi)| \ge \inf_{z\in \partial B(x,2\eta)} |m_n(\xi)| ) = O(\delta^{\kappa \epsilon_0}).$$

Using Cauchy's theorem and arguing as in the proof of Claim~\ref{cl.deviationE}, we obtain 
\begin{eqnarray*}
\P(\sup_{\xi\in \partial B(x,2\eta)} |r_{n,G}(\xi)| \ge \lambda  \sup_{\xi\in B(x,3\eta)} \sqrt{Var[r_{n,G}(\xi)]}) &\lesssim& e^{-\alpha_0 \lambda^2},
\end{eqnarray*}
for some $\alpha_0>0$ and any $\lambda>0$.
Using Lemma~\ref{l.elementary} and Condition~\ref{cond.polyvar}, we also have
$$\sup_{\xi\in B(x,3\eta)} \sqrt{Var[r_{n,G}(\xi)]} \approx \inf_{\xi\in B(x,3\eta)} \sqrt{Var[r_{n,G}(\xi)]}.$$
Thus,  using the given hypothesis we obtain, for some $c'>0$,
\begin{eqnarray*}
&&\P(\sup_{\xi\in \partial B(x,2\eta)} |r_{n,G}(\xi)| \ge t \inf_{\xi\in B(x,2\eta)} |m_n(\xi)|) \\
&\lesssim& \P(\sup_{\xi\in \partial B(x,2\eta)} |r_{n,G}(\xi)| \gtrsim C' t |\log \delta|^{1/2} \sup_{\xi\in B(x,3\eta)} \sqrt{Var[r_{n,G}(\xi)]})\\
&\lesssim& e^{-c'(C't)^2|\log \delta|}.
\end{eqnarray*}
Let $t=1$ in the last estimate. Then  for any $\kappa>0$ we could choose $C_0 \approx \sqrt \kappa$ but large such that this estimate is bounded above by $O((\eta/\delta)^{\kappa})$, as desired.

\subsubsection{Proof of Lemma~\ref{l.nearR-mlarge}, part (ii)} The proof is entirely similar to part (i).

%

 \section{Proof of universality    for real correlation functions}

Below we prove   part (i) of Theorem~\ref{t.corr}, and the same argument may be used to prove part (ii) of this theorem (details will be omitted).

Let $x=(x_1,\dots, x_m) \in I_{\R}(\delta)^m$ and $z=(z_{m+1},\dots, z_{m+k})\in I_{\C_+}(\delta)^k$.  For convenience of notation write $z_j = x_j+ i y_j$ for all $j$. Then for $j \le m$ we have   $y_j=0$ and $x_j\in I(\delta)$, while for $j>m$ we have $y_j>0$. Note that $x_{j}$ and $y_{j}$ may not be inside $I_{\C}(\delta)$ for $j>m$.

Arguing as in the proof of Theorem~\ref{t.complex-corr},   
  it suffices to show that
$$|\E (\prod_{j=1}^{m+k} X_{j})  - \E (\prod_{j=1}^{m+k} X_{G,j})| \lesssim \delta^c, \quad \text{where } \  X_{j} = \begin{cases} \sum_{\alpha \in Z\cap \R} F_{j,\delta}(\alpha-z_j), & j\le m;\\
\sum_{\alpha \in Z \cap \C_+} H_{j,\delta}(\alpha-z_j), & j>m.\end{cases}$$
($X_{G,j}$ are Gaussian analogues), and $F_{j,\delta}$ and $H_{j,\delta}$ satisfy the following conditions:

(i) for each $j\le m$, $F_{j,\delta}$ is in $C^2(\R)$, supported in $(-c\delta, c\delta)$ such that $|F^{(\ell)}_{j,\alpha}|\le 1$ for $\ell=0,1,2$. 

(ii) for each $j>m$, $H_{j,\delta}$ is supported on $B_{\C}(0, c\delta)$ and  is also $C^2(\R^2)$ with $|\partial^{\alpha} H_{j,\delta}| \le \delta^{-|\alpha|}$ for any $|\alpha|\le 2$.

Let $\epsilon_0>0$ be sufficiently small, as required by  Lemma~\ref{l.nearR} and let $\eta=\delta^{1+\epsilon_0}$. 

Let $c'\in (0,1)$ be small such that $c+c'<\widetilde c$.

Let $\Phi:\R\to \R$ be a bump function supported on $[-c', c']$ with $\Phi(0)=1$. 

Let $\Psi:\R\to [0,1]$ be  a smooth function supported on $\{x\ge c'/2\}$ such that $\Psi(x)=1$ if $x \ge c'$. 

Let $L=O(1)$ be sufficiently large. Let $K_{1,\delta},\dots, K_{m+k,\delta}:\C\to\C$ be defined by 
$$K_{j,\delta}(x+iy)=\begin{cases}\delta^{L\epsilon_0} F_{j,\delta}(x) \Phi(y/\eta), & j\le m,\\ 
 \delta^{L\epsilon_0}H_{j,\delta}(x+iy)\Psi( (y + y_{j})/\eta), & j\ge m+1.\end{cases}$$ One could check that $K_{1,\delta},\dots, K_{m+k,\delta}$ are supported on $B(0,(c+c')\delta)$ and are   $C^2(\R^2)$ with  $\partial^\alpha$   derivatives  bounded by $O(\delta^{-|\alpha|})$ for any multi-index $|\alpha|\le 2$.

Applying Theorem~\ref{t.complex-corr} for test functions of tensor-product type, it follows that for some $\alpha_0>0$ (which does not depend on  $\epsilon_0$) we have
$$|\E(\prod_{j=1}^{m+k} Y_{j}) - \E (\prod_{j=1}^{m+k} Y_{G,j})| \lesssim \delta^{\alpha_0}, \quad \text{where }   \ Y_{j}(z) := \sum_{z\in Z} K_{j,\delta}(z-z_j).
$$
Letting $Z_{j} := \delta^{-L\epsilon_0}Y_{j}$ and making sure  $\epsilon_0<c_0/(Lm+Lk)$, 
 it remains to show
$$\E |\prod_{j=1}^{m+k} X_{j}  - \prod_{i=1}^{m+k} Z_{j}| = O(\delta^{\alpha_1})$$
for some $\alpha_1>0$.
Since $X_{j}, Z_j  \le N_{p_n}(B(z_j, c\delta))$, using Corollary~\ref{c.localcount} we have  
$\E |X_{j}|^{m+k},  \E |Z_{j}|^{m+k}  \lesssim  |\log \delta|^{O(m+k)}$. 
Via  Holder's inequality, it therefore suffices to show that for  some $c>0$ we have 
$$\E |X_{j} - Z_{j}|^{m+k} \lesssim \delta^{c}.$$
Now, for each $1\le j\le m+k$ let 
\begin{align*}
&&S_j  =  \{t\in \R: \ |t- sign(Re(z_j))|z_j||\le (c+c')\delta\}  \times [-c'\eta, c'\eta].
\end{align*}
We first show that if $X_j-Z_j \ne 0$ then $|Im(z_j)| \le (c+c')\delta$ and
\begin{eqnarray}\label{e.Xj-Zj}
|X_j -Z_j|  &\lesssim& |Z \cap (S_j \setminus \R)|.
\end{eqnarray}

Indeed, we first consider $1\le j\le m$. Then $z_j=x_j\in I_{\R}(\delta)$. Therefore,
\begin{eqnarray*}
X_{j} - Z_{j} &=& \sum_{\alpha \in Z\cap \R} F_{j,\delta}(\alpha-x_j) - \sum_{\alpha\in Z} F_{j,\delta}(Re(\alpha)-x_j) \Phi(Im(\alpha) / \eta)\\
&=& -\sum_{\alpha\in Z\setminus \R}   F_{j,\delta}(Re(\alpha)-x_j)  \Phi(Im(\alpha) / \eta) \quad \text{(since $\Phi(0)=1$)}. 
\end{eqnarray*}
Since both $F_{j,\delta}$ and $F_{j,\delta}$ are bounded, it suffices to show that any $\alpha$ that contributes to the sum must be in $S_j$. Indeed, for such $\alpha$ we  have $|Re(\alpha)-x_j| < c\delta$ and  $|Im(\alpha)| < c'\delta$, which implies the desired claim.

We now consider $m+1\le j\le m+k$.  We have
$$X_{j} - Z_{j}  = \sum_{\alpha\in Z\cap \C_+} H_{j,\delta}(\alpha-z_j)  - \sum_{\alpha \in Z} H_{j,\delta}(\alpha-z_j) \Psi(Im(\alpha)/\eta).$$
Since $\Psi$ is supported on $[c'/2,\infty)$  in the second summation we could further assume that $\alpha \in \C_+$. We obtain
$$X_{j} - Z_{j}  = \sum_{\alpha\in Z\cap \C_+} H_{j,\delta}(\alpha-z_j) (1- \Psi(Im(\alpha/\eta)) ).$$
For any contributing $\alpha$,   it holds that $|Im(\alpha)| < c'\eta$, therefore 
$$|Im(z_j)|\le |Im(\alpha)|+|Im(\alpha)-Im(z_j)| <  (c+c')\delta.$$ 
In particular, $|Re(z_j)| \ge |z_j| - |Im(z_j)| \ge 1- O(\delta)$ and this can be made very large compared to $\delta$. Now, 
$$|Re(\alpha)-Re(z_j)|\le |\alpha-z_j|\le c\delta$$
therefore  $Re(\alpha)$ has the same sign as $Re(z_j)$. Thus it remains to show that $||Re(\alpha)|-|z_j|| \le (c+c')\delta$. Now, using the triangle inequality this follows from 
$$||Re(\alpha)|-|z_j|| \le  ||\alpha|-|z_j||+ |Im(\alpha)| \le |\alpha-z_j|+|Im(\alpha)| \le (c+c')\delta.$$
This completes the proof of \eqref{e.Xj-Zj}.

Now, the strip $S_j$ could be covered by $O(\delta^{-\epsilon_0})$ sets of the form $B(x, \eta)$ with center $x$ inside  $(sign(Re(z_j))|z_j|-(c+c')\delta, sign(Re(z_j))|z_j|+(c+c')\delta)$. Since $c+c'<\widetilde c$ and since $Im(z_j)|\le (c+c')\delta$, it follows that for such $x$ the ball $B(x,2\eta)$ would be inside the interval $J$ where the given hypothesis on the relationship between $m_n$ and $r_n$ holds. Now, since $p_n$ is a real polynomials its complex roots are symmetric about the real axis. Thus, using the small ball estimates proved in Lemma~\ref{l.nearR} (if $m_n$ is small compared to $r_n$) or the small ball estimates proved in Lemma~\ref{l.nearR-mlarge} (if $m_n$ is large compared to $r_n$) with  $\kappa=3/2$, together with an union bound, we obtain
$$P(|Z\cap (S_j\setminus \R)|\ge 1) = O(\delta^{-\epsilon_0}\delta^{3\epsilon_0/2}) = O(\delta^{\epsilon_0/2}).$$

Now, since $|Z\cap (S_j\setminus \R)|$ is a nonnegative integer, by Theorem~\ref{t.localcount} we have
\begin{eqnarray*}
\E |X_{j} - Z_{j}|^{m+k} 
&\lesssim& \E [1_{ |Z\cap (S_j\setminus \R)|\ge 1} N_{p_n}(B(z_j, c\delta))^{m+k}] \\
&\lesssim& \delta^{\epsilon_0/2} |\log \delta|^{O(m+k)} \quad \lesssim \quad \delta^{\epsilon_0/3}.
\end{eqnarray*}
This completes the proof of Theorem~\ref{t.corr}.

\endproof

 \section{Reduction of Theorem~\ref{t.general}  to Gaussian polynomials} 
In this section,  using Theorem~\ref{t.corr}  we will reduce   Theorem~\ref{t.general} to Gaussian random polynomials. The proof of   Theorem~\ref{t.general} for Gaussian polynomials will be discussed in the next section.

Let $B_C=\{1-\frac 1 C \le |t| \le 1+\frac 1C\}$. Using Lemma~\ref{l.im}, to reduce Theorem~\ref{t.general} to the Gaussian setting, it suffices to show that
\begin{eqnarray*}
|\E N_n(I\cap B_C) - \E N_{G,n}(I\cap B_C)| &=& O(1) .
\end{eqnarray*}
Thus without loss of generality we may assume that   $I\subset [1-1/C, 1+1/C]$ or $I\subset [-1-1/C, -1+1/C]$. Below, we will only consider the first case, and we may use the same argument for the other case.

Let $\epsilon>0$ be a very small absolute constant. Recall the definition of $I(\delta)$   from \eqref{e.Idelta} and the paragraph after \eqref{e.Idelta}. Let $\widetilde I_{\R}(\delta) = \{z: 1/z\in I_{\R}(\delta)\}$.

Note that  we may cover  $I$ using intervals    $I_{\R}(2^m)$ and $\widetilde I_{\R}(2^\ell)$ where $\frac 1 {n}\lesssim 2^{m} \lesssim \frac 1 C$ and $\frac 1 n \lesssim 2^\ell \lesssim \frac 1 C$. Let $M$ and $L$ be respectively the sets of $m$ and $\ell$ such that $I(2^m)$ and $\widetilde I(2^{\ell})$  intersect $I$.  Clearly,  nearby covering intervals have comparable lengths. Thus, we may construct a sequence of functions $\varphi_m, \psi_{\ell}$ (similar to a partition of unity) such that  $\varphi_m$ is supported on $(1+\epsilon) I(2^{m})$  and $\psi_\ell$ is supported on $(1+\epsilon) \widetilde I(2^{\ell})$, and furthermore

(i) $|\partial^{\alpha}\psi_\ell| \lesssim 2^{|\alpha|\ell}$ and $|\partial^{\alpha} \varphi_m|\lesssim 2^{|\alpha|m}$ for any partial derivatives, and 

(ii) $\gamma(y):=\sum_{m\in M}  \varphi_m(y) + \sum_{\ell\in L} \psi_{\ell}(y)$
is equal to $1$ for all $y\in I$ and is supported inside $I\cup I_l\cup I_r$ where $I_l, I_r$ are two intervals from the covering  that contain endpoints of $I$.  

Now, we could shrink the endpoint intervals   $I_l$ and $I_r$ by factors comparable to $1$ (if necessary) so that $I$ remains covered by the new collection of intervals, and at the same time $(1+2\epsilon)I_l, (1+2\epsilon)I_r$ are subsets of the assumed enlargement $J$ of $I$. The given definition of enlargement ensures that  the shrinking of these intervals could be done. We may redesign the bump functions $\phi_m$ and $\psi_\ell$  associated with $I_l$ and $I_r$ such that they will still be supported inside $(1+\epsilon)I_l$ and $(1+\epsilon)I_r$, respectively.

It follows from Theorem~\ref{t.corr} that, for some $\alpha_1>0$,
$$|\E \sum_{\alpha\in Z\cap \R} \varphi_m(\alpha) - \E\sum_{\alpha\in Z_G\cap \R}\varphi_m(\alpha) | = |\int_{\R} \varphi_m(y)[d\sigma(y) - d\sigma_G(y)]| \lesssim  2^{m \alpha_1},$$
$$|\E \sum_{\alpha\in Z\cap \R} \psi_{\ell}(\alpha) - \E\sum_{\alpha\in Z_G\cap \R}\psi_{\ell}(\alpha) | = |\int_{\R} \psi_{\ell}(y)[d\sigma(y) - d\sigma_G(y)]| \lesssim  2^{\ell \alpha_1}.$$

Summing the last two estimates over $m$ and $\ell$, we obtain
$$|\E \sum_{\alpha\in Z\cap \R} \gamma(\alpha) - \E\sum_{\alpha\in Z_G\cap \R}\gamma(\alpha) | = O(1).$$

Now, $|\E N_n(I) - \E \sum_{\alpha\in Z\cap \R} \gamma(\alpha)| = O(\E N_n(I_l \cup I_r))$. For the local intervals $I_l$ and $I_r$,  we will show that $\E N_n(I_l)=O(1)$ and $\E N_n(I_r)=O(1)$. Since the details are entirely similar we will only discuss the estimate for $\E N_{n}(I_l)$. Since $(1+\epsilon)I_l\subset J$ the enlargement of $I$, we may construct a bump function $\phi$ adapted to $I_l$ that equals $1$ on $I_l$ but vanishes outside $(1+  \epsilon / 2)I_l$, in particular its support is strictly contained inside $J$. Let $d\rho$ be the $1$-point correlation measure for the real root of $p_n$ and $d\rho_G$ be its Gaussian analogue. By Theorem~\ref{t.corr}, we obtain
$$\E N_{n}(I_l) \le \int \phi d\rho = \int \phi d\rho_G + O(1) \le \E N_{n,G}((1+\epsilon/2)I_l) + O(1)$$
Then assuming that the Gaussian case of Theorem~\ref{t.general} is known and using the fact that $J$ remains an enlargement of $(1+\epsilon/2)I_l$, we obtain
$$|\E N_{n,G}((1+\epsilon/2)I_l)|  \le |\E N_{r_{n,G}} ((1+\epsilon/2)I_l)| + O(1) = O(1)$$
here in the last estimate we may use  Proposition~\ref{p.S} in the next section (which is a consequence of explicit Gaussian computations in \cite{dnv2017}). 

This completes the proof of the reduction of Theorem~\ref{t.general} to  Gaussian polynomials.

\section{Proof of Theorem~\ref{t.general}  for Gaussian polynomials} \label{s.gaussian}
In this section we prove Theorem~\ref{t.general} for the Gaussian  polynomial $p_n(t)=\sum_{j=0}^n (b_j+c_j\xi_j)t^j$ where $\xi_j$ are iid normalized Gaussian, and throughout the section we will assume that $b_j$ and $c_j$ satisfy Condition~\ref{cond.polyvar}. 

Let $m_n =\E [p_n]$ and $r_n(t)=\sum_j c_j \xi_j t^j$ and let  $\mathcal P = Var[r_n(t)]$, $\mathcal Q  = Var[r'_n(t)]$, and $\mathcal R = Cov[r_n(t),r'_n(t)]$, and $\mathcal S  = \mathcal P \mathcal Q - \mathcal R^2$. 


We recall the following Kac-Rice formula \cite[Corollary 2.1]{farahmand1998}. Let $erf(x) =\int_0^x e^{-t^2}dt$. Then $\E N_n(a,b) = I_1(a,b) + I_2(a,b)$ where

\begin{eqnarray} 
\label{e.I1def} I_1(a,b) &=& \int_a^b \frac{\mathcal S^{1/2}}{\pi \mathcal P} \exp(-\frac{m_n^2 \mathcal Q + m_n'^2 \mathcal P - 2m_n m_n' \mathcal R}{2\mathcal S})dt \\
\label{e.I2def} I_2(a,b) &=& \sqrt 2 \int_a^b \frac{|m_n'\mathcal P - m_n \mathcal R|}{\pi \mathcal P^{3/2}} \exp(-\frac{m_n^2}{2\mathcal P}) erf(\frac{|m_n'\mathcal P - m_n \mathcal R|}{\sqrt {2\mathcal P}\mathcal S}) dt. 
\end{eqnarray}

We will also work with the normalized reciprocal polynomial $  p^*_n(t)= m^*_n(t) + r^*_n(t)$, and we will denote by $I_1^*$, $I_2^*$, $\mathcal P^*, \mathcal Q^*,\mathcal R^*,\mathcal S^*$ the analogous quantities.  
 
Using Lemma~\ref{l.im},  we may assume without loss of generality that $I \subset \{1-c\le |t|\le 1+c\}$ for a (small) absolute constant $c>0$. By breaking up $I$ into $I_{>1}$ and $I_{\le 1}$ and notice that $N_{p_n}(I_{>1}) = N_{p^*_n}(K)$ where $K=\{1/t, \ \ t\in I_{>1}\}$ we may reduce the consideration to $I\subset \{1-c\le |t|\le 1\}$.

Now, using Lemma~\ref{l.elementary}, we have

\begin{corollary} \label{c.PQR} Assume that $b_j$ and $c_j$ satisfy Condition~\ref{cond.polyvar}. Then for any $c\in(0,1)$ it holds uniformly over $1-c\le |t|\le 1$ that
\begin{eqnarray*}
\mathcal P(t)   \approx (1+1/n-|t|)^{-(2\rho+1)}, \quad \mathcal P^*(t) \approx (n+1)^{2\rho} (1+1/n-|t|)^{-1}, \\
\mathcal Q(t)  \approx (1+1/n-|t|)^{-(2\rho+3)}, \quad \mathcal Q^*(t) \approx (n+1)^{2\rho}(1+1/n-|t|)^{-3}, \\
|\mathcal R(t)|  \approx (1+1/n-|t|)^{-(2\rho+2)}, \quad \mathcal R^*(t) \approx (n+1)^{2\rho}(1+1/n-|t|)^{-2}.
\end{eqnarray*}
\end{corollary}

On the other hand, by the classical Kac formula,   $\rho_n(t):=\frac{\mathcal S^{1/2}}{\pi\mathcal P}$ is the density for the real root distribution of  $r_n(t)=\sum_j c_j \xi_j t^j$, and similarly $\rho^*_n(t):=\frac{{\mathcal S^*}^{1/2}}{\pi {\mathcal P}^*}$ is the density for the real root distribution of  $r^*_n(t)$, and both of them can be easily bounded by $O(n)$ by elementary inspection.  Note that the Gaussian density for $r_n(t)=\sum_j c_j \xi_j t^j$ (and for its reciprocal polynomial) was studied\footnote{In fact, in \cite{dnv2017}  it was required that $|c_j|\sim (1+j)^{\rho}$ for $O(1)\le j\le n$,  however the Gaussian computations in \cite{dnv2017} can be easily modified to work with the weaker assumption $O(1)\le j\le n-O(1)$ in the current paper.} in \cite{dnv2017}, and we summarize the known estimates for them from \cite[Lemma 10.3, Lemma 10.6]{dnv2017} in the following proposition. 
\begin{proposition}\label{p.S}
Assume that  $c_j$ satisfy Condition~\ref{cond.polyvar}.  Let $c>0$ be small. Then uniformly over $1-c\le |t|\le 1-c'/n$ we have
\begin{eqnarray*}
\mathcal S(t)   \approx \mathcal P(t)^2 (1-|t|)^{-2}, \quad \mathcal S^*(t) \approx \mathcal P^*(t)(1-|t|)^{-2},
\end{eqnarray*}
and uniformly over $1-c'/n \le |t|\le 1+c'/n$ we have
\begin{eqnarray*}
\mathcal S(t)   \lesssim n^2 \mathcal P(t)^2 , \quad \mathcal S^*(t)  \lesssim n^2 \mathcal P^*(t) .
\end{eqnarray*}
\end{proposition}
In fact, in the original setting considered in \cite{dnv2017} it was  required that $c_j\approx (1+j)^\rho$ for all $O(1)\le j\le n$, so it is a little stricter than our setting $O(1)\le j\le n-O(1)$, however  the computation in the Gaussian setting in \cite{dnv2017} is not affected much with our slightly more relaxed assumption. We omit the details.

\subsection{Estimates for $I_2$}  
We will show that, under the hypothesis of Theorem~\ref{t.general} about the relative relation between $m_n$ and $r_n$ on $I$, we will always have $I_2(I)=O(1)$. We separate the proof into two cases, depending on whether $m_n$ dominates $r_n$ or is dominated by $r_n$.

First, we consider the situation when the deterministic component $m_n$ dominates the random component $r_n$ on $I$.

\begin{lemma} \label{l.basicI2-large}  Let $c>0$. There is a constant $C>0$ such that the following holds.  Let $I \subset \{1-c\le |t|\le 1\}$  be an interval whose endpoints may depend on $n$. 
\begin{align*}
\text{(i) Assume that}&&
 |m_n(t)| \ge C |\log (1+\frac 1 n-|t|)|^{1/2} \sqrt{Var[r_n(t)]} \quad \text{ for $t\in I$}.
 \end{align*}
\begin{align*}
\text{Then} && I_2(I)  = O(1).
\end{align*}
\begin{align*}
\text{(ii) Assume that}&&  |m^*_n(t)| \ge C |\log (1+\frac 1 n-|t|)|^{1/2} \sqrt{Var[r^*_n(t)]} \quad \text{for $t\in I$}.
\end{align*}
\begin{align*}
\text{Then} && I^*_2(I)  = O(1).
\end{align*}

\end{lemma}

\proof  Using Lemma~\ref{l.elementary} and Corollary~\ref{c.PQR} we have 
$$\frac{|m_n'\mathcal P - m_n \mathcal R|}{\mathcal P^{3/2}} \lesssim \frac{|m_n'|}{\mathcal P^{1/2}} + \frac{|m_n \mathcal R|}{\mathcal P^{3/2}}  \lesssim (1+\frac 1 n-|t|)^{-3/2} $$
and by the given hypothesis  $|m_n(t)|^2/\mathcal P \ge 2 C' |\log (1+\frac 1 n -|t|)|$ where $C'$ is comparable to $C^2$. Therefore
\begin{eqnarray*}
I_2(I) &\lesssim&   \int_I (1-|t|+\frac 1 n)^{C'-3/2}  dt 
\end{eqnarray*}
so  if $C$ is big enough then $C'>5/2$ and the last   integral is $O(1)$, as desired.

The consideration for $I^*_2(I)$ is entirely similar.
\endproof

We now consider the situation when  $m_n$ is dominated by  $r_n$. 

%

Recall that  $\phi:(0,1)\to [0,1]$ is such that the following holds for some $c>0$:
\begin{eqnarray}\label{e.phi}\int_{1/n}^{c} \frac {\phi(t)}{t}dt &=& O(1).
\end{eqnarray}

\begin{lemma} \label{l.basicI2-small} Let $c>0$ and let  $\phi:(0,1)\to \R_+$ satisfy \eqref{e.phi}. Let $I \subset \{1-c\le |t|\le 1\}$  be an interval whose endpoints may depend on $n$. 

(i)  Assume that the following holds uniformly over $t\in I$.
$$|m_n(t)| \lesssim \phi(1-|t|+\frac 1 n) \sqrt{Var[r_n(t)]}, $$
$$|m'_n(t)| \lesssim  \phi(1-|t|+\frac 1 n)  \sqrt{Var[r'_n(t)]}.$$
\begin{align*}
\text{Then}&& I_2(I)=O(1).
\end{align*}
(ii) Under the analogous assumptions, we also have $I^*_2(I)  = O_\epsilon(1)$.

\end{lemma}

\proof Using the given hypothesis and using Corollary~\ref{c.PQR}, we have 
$$\frac{|m_n'\mathcal P - m_n \mathcal R|}{\mathcal P^{3/2}}  \ \lesssim  \   \phi(1-|t|+\frac 1n)(\frac{\mathcal Q^{1/2}}{\mathcal P^{1/2}} + \frac{ \mathcal R}{\mathcal P}) \  \lesssim \  \frac   {\phi(1-|t|+\frac 1n)}{1-|t|+\frac 1n}.$$
Since $\exp(-m_n^2/\mathcal P)\le 1$, we obtain 
$$I_2(I) \lesssim \int_{1-c}^{1} \frac   {\phi(1-t+\frac 1n)}{1-t+\frac 1n} dt \le \int_{1/n}^{c+1/n} \frac{\phi(t)}{t}dt = O(1).$$
This completes the proof of part (i). The second part (ii) can be proved similarly.
\endproof

\subsection{Estimates for $I_1$}

Here we will also divide the consideration into two cases, depending on whether $m_n$ is dominant or $r_n$ is dominant.

The following result addresses   the situation when $m_n$ is dominated by $r_n$. 

\begin{lemma}\label{l.basicI1-small} Assume that $\phi: (0,1)\to\R_+$ satisfies \eqref{e.phi}.  Let $c>0$ and let $I\subset \{1-c\le |t|\le 1\}$ be an interval whose endpoints may depend on $n$. 

(i) Assume that uniformly over $t\in I$ we have
$$|m_n(t)| \le \sqrt{\phi(1-|t|+\frac 1n)} \sqrt{Var[r_n(t)]},$$  
$$|m'_n(t)| \le  \sqrt{\phi(1-|t|+\frac 1 n)} \sqrt{Var[r'_n(t)]}$$
\begin{align*}
\text{Then} &&  I_1(I)  = \int_I \rho_n(t) dt + O(1).
\end{align*}

(ii) Under analogous assumptions, a similar estimate holds for $I^*_1(I)$.

\end{lemma}

The following result deals with the situation when $m_n$ dominates $r_n$.

\begin{lemma}
 \label{l.basicI1-large}  Let $c>0$ and let $I\subset \{1-c\le |t|\le 1\}$ be an interval whose endpoints may depend on $n$. 

(i) Assume that uniformly over $t\in I$ we have
\begin{eqnarray*}
|m_n(t)| &\gtrsim&   |\log (1-|t|+\frac 1 n)|^{1/2} \sqrt{Var[r_n(t)]}.
\end{eqnarray*} 
\begin{align*}
\text{Then} && I_1(I)=O(1).
\end{align*}

(ii) Under analogous assumptions, a similar estimate holds for $I^*_1(I)$.
\end{lemma}

The proof of these results are based on the following technical estimate.
For convenience, let $\mathcal T(t)=\frac{m_n^2}{\mathcal P} + \frac{m_n'^2}{\mathcal Q}$, and define $\mathcal T^*(t)$ analogously. Recall that $ \rho_n(t):=\frac{\mathcal S^{1/2}}{\pi\mathcal P}$ is the density for the real root distribution of  $r_n(t)$, and $\rho^*_n:=\frac{{\mathcal S^*}^{1/2}}{\pi\mathcal P^*}$ is the density for the real root distribution for $r^*_n$.

\begin{lemma}\label{l.basicI1}  Let $c>0$ be sufficiently small and let $c'>0$ be sufficiently large. Then there are finite absolute constants $C_1, C_2>0$ that may depend  on $c,c'$ such that the following holds for any interval $I$ whose endpoints may depend on $n$. 

(i) If $I\subset \{1-c'/n \le |t|\le 1+c'/n\}$ then $I_1(I)=O(1)$ and $I^*_1(I)=O(1)$.

(ii) If $I\subset  \{1-c\le |t|\le1-c'/n\}$ then 
\begin{align*}
&& \int_I   \rho_n(t) e^{-C_1 \mathcal T(t)} dt \quad \le   \quad I_1(I) \quad \le   \quad \int_I    \rho_n(t) e^{-C_2 \mathcal T(t)} dt,
 \end{align*}
 and the analogous estimate holds for $I^*_1(I)$.

\end{lemma}

\proof  (i) Since $\mathcal P\mathcal Q \ge \mathcal R^2$, it follows that  $m_n^2\mathcal Q + m_n'^2 \mathcal P - 2 m_n m_n'\mathcal R\ge 0$, so 
$$I_1(I)  \lesssim \int_{||t|-1|\lesssim 1/n} \rho_n(t) dt   = O(1).$$
The estimate for $I^*_1$ is proved similarly.

(ii) Let $1-c\le |t|\le 1-c'/n$. From Corollary~\ref{c.PQR} and Proposition~\ref{p.S}, we obtain
\begin{eqnarray*}
\mathcal P^{1/2}\mathcal Q^{1/2}- |\mathcal R|
&=& \frac{\mathcal S}{\mathcal P^{1/2}\mathcal Q^{1/2}+|\mathcal R|}  
\ \gtrsim \  (1-|t|)^{-(2d+2)} \gtrsim |\mathcal R|.
\end{eqnarray*}
In other words for some $C>0$ we have $\mathcal P^{1/2}\mathcal Q^{1/2} \ge (1+C)|\mathcal R|$. Consequently, by the geometric mean inequality we have 
\begin{eqnarray*}
 m_n^2\mathcal Q + m_n'^2 \mathcal P - 2m_n m_n'\mathcal R \approx  m_n^2\mathcal Q + m_n'^2 \mathcal P.
\end{eqnarray*}
Now, by Corollary~\ref{c.PQR} and Proposition~\ref{p.S} we have $\mathcal S \approx \mathcal P \mathcal Q$. It follows that
\begin{eqnarray*}
\frac{m_n^2\mathcal Q + m_n'^2 \mathcal P - 2m_n m_n'\mathcal R}{\mathcal S} \ \ \approx \ \  \frac{m_n^2\mathcal Q + m_n'^2\mathcal P}{\mathcal P \mathcal Q} &=&  \mathcal T(t).
\end{eqnarray*}
The desired estimate then follows from the definition \eqref{e.I1def} of $I_1$.

The proof for $I^*_1(t)$ is completely analogous.
\endproof
 We now use Lemma~\ref{l.basicI1} to prove Lemma~\ref{l.basicI1-small} and Lemma~\ref{l.basicI1-large}. Below we will show only the proof for the desired estimates for $I_1$, the same argument works for $I^*_1$. We start with the case when $m_n$ is dominated by $r_n$: under the assumptions of  Lemma~\ref{l.basicI1-small} we have  $\mathcal T(t) \lesssim \phi(1-|t|+\frac 1 n)$. Using $1\ge e^{-x}\ge 1-x$ for $x\ge 0$  and using Proposition~\ref{p.S}, it follows that
 $$|I_1(I) -   \int_I \rho_n(t)dt| \lesssim |\int_I \rho_n(t) \mathcal T(t)dt|  \lesssim  \int_{1-c}^1 \frac{\phi(1-t+\frac 1 n)}{1-t+\frac 1 n}dt =   O(1).$$
Now in the case when $m_n$ dominates $r_n$: under the assumptions of Lemma~\ref{l.basicI1-large}  we have $\mathcal T(t) \gtrsim |\log(1+\frac 1 n -|t)|$, while $\rho_n(t)\lesssim (1+\frac 1 n - |t|)^{-1}$ thanks to Proposition~\ref{p.S}. Therefore, for some $c''>0$ we have
$$I_1(t) \lesssim \int_{1-c}^{1} \frac 1 {1-t+\frac 1 n} e^{-c'' |\log (1-t+\frac 1 n)|}dt  = \int_{1/n}^{c+1/n} u^{c''-1}du = O(1).$$

\appendix

\end{document}